\documentclass[11pt,draft]{amsart}  
\usepackage{amssymb,amsxtra} 
\usepackage{verbatim,fullpage,enumerate} 
\usepackage[french,english]{babel}
\usepackage[T1]{fontenc}
\usepackage{xspace}

\numberwithin{equation}{section}

\newtheorem{Theoreme}{Th{\'e}or{\`e}me}[section] 
\newtheorem{Prop}[Theoreme]{Proposition}
\newtheorem{prop}[Theoreme]{Proposition} 
\newtheorem{Assertion}[Theoreme]{Assertion} 
\newtheorem{Lemme}[Theoreme]{Lemme} 
\newtheorem{lem}[Theoreme]{Lemme}
\newtheorem{cor}[Theoreme]{Corollaire}

\theoremstyle{definition}

\theoremstyle{remark}
\newtheorem{Remarque}[Theoreme]{Remarque}

\newenvironment{dem}{\proof[D{\'e}monstration]}{\endproof}

\def\preisomto{\vbox{\hbox to 
                 15pt{\hfill$\sim$\hfill}\nointerlineskip\vskip 
                 -0.3pt  
                 \hbox to 15pt{\rightarrowfill}}} 
 
\def\isomto{\mathop{\preisomto}} 
 
\def\prelongisomto{\vbox{\hbox to
                18pt{\hfill$\sim$\hfill}\nointerlineskip\vskip -0.3pt
                \hbox to 18pt{\rightarrowfill}}}  

\def\longisomto{\mathop{\prelongisomto}}

\def\Ker{\operatorname {Ker}}

\def\Hom{\operatorname {Hom}}
\def\End{\operatorname {End}}
\def\Aut{\operatorname {Aut}}

\def\Int{\operatorname {Int}}
\def\id{\operatorname {id}}  
\def\lg{\operatorname {lg}}
\def\sgn{\operatorname {sgn}}
\def\GL{\operatorname {GL}}

\newcommand{\sto}{\twoheadrightarrow}

\newcommand{\N}{{\mathbb{N}}}
\newcommand{\C}{{\mathbb{C}}}
\newcommand{\R}{{\mathbb{R}}}

\newcommand{\D}{{\mathbb{D}}}
\renewcommand{\H}{{\mathbb{H}}}

\newcommand{\hb}{\mathbb{H}_{B}}
\newcommand{\hd}{\mathbb{H}_{D}}

\newcommand{\hbi}{\overline{\mathbb{H}_B(\mu_i)}}
  
\newcommand{\VC}{{V_{\mathbb{C}}}}
\newcommand{\CW}{{\mathbb{C} W}}
\newcommand{\CT}{{\mathbb{C} T}}
\newcommand{\CU}{{\mathbb{C} U}} 

\newcommand{\aalg}{\mathfrak{a}}
\newcommand{\fa}{\mathfrak{a}}

\newcommand{\che}{\mathcal{H}}
\newcommand{\calP}{\mathcal{P}}
\newcommand{\calS}{\mathcal{S}}

\newcommand{\calC}{\mathcal{C}}

\newcommand{\calA}{\mathcal{A}}

\newcommand{\calH}{\mathcal{H}}

\newcommand{\calL}{\mathcal{L}}

\newcommand{\mbH}{\mathbf{H}}

\newcommand{\DK}{\hat{D}} \newcommand{\tdel}{\tilde{\Delta}}
\newcommand{\tk}{\tilde{k}} \newcommand{\bk}{\bar{k}}
\newcommand{\tgamma}{{\tilde{\gamma}}}
\newcommand{\bgamma}{{\bar{\gamma}}}

\newcommand{\tnu}{{\tilde{\nu}}}
\newcommand{\bmu}{\bar{\mu}} 
\newcommand{\bnu}{\bar{\nu}}  
\newcommand{\vphi}{\varphi}
\newcommand{\vt}{\vartheta}
\newcommand{\vkappa}{\varkappa}

\newcommand{\half}{{\frac{1}{2}}}

\newcommand{\twisted}[2]{{}^{{#1}}\!{{#2}}} 
\newcommand{\RMod}[1]{{{#1}\text{-}\mathrm{Mod}}}
\newcommand{\Rmod}[1]{{{#1}\text{-}\mathrm{mod}}}
\newcommand{\ModR}[1]{{\mathrm{Mod}\text{-}{#1}}}
\newcommand{\modR}[1]{{\mathrm{mod}\text{-}{#1}}}

\newcommand{\dual}[2]{{\langle {#1} ,  {#2} \rangle}}

\begin{document}

{\selectlanguage{french}

  \title[]{Une g\'en\'eralisation de l'alg{\`e}bre de Hecke
    gradu\'ee de type~$\mathbf{B}$}
  \author[C.~Dez{\'e}l{\'e}e]{Charlotte Dez{\'e}l{\'e}e}
  \address{D{\'e}partement de math{\'e}matiques, Universit{\'e}
    de Brest, 29285 Brest cedex, France}
  \email{Charlotte.Dezelee@univ-brest.fr} \date{11 juin 2003}

\begin{abstract} 
  On \'etudie une alg\`ebre ``proche'' d'une alg\`ebre de Hecke
  gradu\'ee et r\'ealisable dans l'alg\`ebre de Cherednik
  rationnelle de type $B$.  On introduit pour cette alg\`ebre
  des modules de la s\'erie principale et l'on d\'emontre un
  crit\`ere d'irr\'eductibilit\'e de ces modules. On en déduit
  des résultats similaires pour une algèbre associée à un
  système de racines de type $D$.
  
  \medskip
  
  For a root system of type $B$ we study an algebra similar to
  a graded Hecke algebra, isomorphic to a subalgebra of the
  rational Cherednik algebra. We introduce principal series
  modules over it and prove an irreducibility criterion for
  these modules. We deduce similar results for an algebra
  associated to a root system of type $D$.
\end{abstract}

\maketitle

\addtocounter{section}{-1}

\section{Introduction}
\label{sec0}

Soit $W \to \GL(\fa_\R^*)$ la représentation << naturelle >>
d'un groupe de Weyl $W$ et $R \subset \fa_\R^*$ le système de
racines associé.  On peut attacher à cette donnée plusieurs
familles de $\C$-algèbres. Nous nous intéressons ici à deux
d'entre elles:
\begin{itemize}
\item L'algèbre de Cherednik rationnelle $\che(k)$,
  (cf.~\cite{TL, EG} pour une définition précise). Elle dépend
  d'une fonction de multiplicités $W$-invariante $k : R \to \C$
  et elle est engendrée par $\fa^* = \C \otimes_\R \fa_\R^*$,
  $W$ et des opérateurs de Dunkl $T_y(k)$ pour $y \in \fa$ (le
  dual de $\fa^*)$.
\item L'algèbre de Hecke graduée $\H(k)$; cette algèbre dépend
  elle aussi d'une multiplicité $k$ et est engendrée par les
  éléments de $\fa$ et $W$, soumis à certaines relations de
  commutation, voir \cite{KR, Lu, Op2}.
\end{itemize}
Les représentations de dimension finie de $\H(k)$ ont été
étudiées dans, par exemple, \cite{KR, Op2} et celles de
$\che(k)$ dans \cite{BEG, TL}.

Il a été remarqué par plusieurs auteurs, e.g.~\cite{ EG, Ka},
que lorsque $R$ est de type $A_{n-1}$ on peut réaliser $\H(k)$
dans $\che(k)$ de la façon suivante. On fixe une base
(orthonormée) $\{e_1,\dots,e_n\}$ de $\fa_\R$, de base duale
$\{z_1,\dots,z_n\}$; alors la sous-algèbre de $\che(k)$
engendrée par $W$ et les $z_i T_{e_i}(k)$, $1 \le i \le n$, est
isomorphe à $\H(k)$. Cette algèbre, que nous noterons $\D$,
existant pour tous les types de systèmes de racines il est
naturel de l'étudier et d'examiner les liens qui peuvent
exister avec les algèbres de Hecke graduées (ou leurs
généralisations).

On se propose ici d'\'etudier les repr\'esentations de
l'algèbre $\D$ pour un système de racines de type $B_n$, dont
les occurrences sont d\'ej\`a nombreuses. S.~Kakei a par
exemple introduit une famille commutative $D_1,\dots,D_n$
d'éléments de $\che(k)$ qui avec $W$ engendrent $\D$. Dans la
première partie nous introduisons (par générateurs et
relations) une algèbre $\H_B(\tk)$ dont nous calculons le
centre et montrons qu'elle est isomorphe à $\D$. Il faut noter
qu'ici la fonction $\tk$ définie sur $R$ est à valeurs dans
l'algèbre du sous-groupe $T$ de $W$ engendré par les réflexions
associées aux racines courtes.

L'étude des représentations de dimension finie des algèbres de
Hecke graduées repose en grande partie sur celle des modules de
la série principale, notés $M(\lambda)$. Ils jouent le rôle des
modules de Verma pour une algèbre de Lie semi-simple et ils
sont paramétrés par les formes linéaires $\lambda \in \fa^*$.
Ainsi, tout $\H(k)$-module irréductible est quotient d'un
$M(\lambda)$. Nous montrons que pour $\H_B(\tk)$ une étude
similaire est possible: on définit au \S~\ref{sec2} des modules
$M(\gamma \otimes \mu)$ indexés par les caractères $\gamma
\otimes \mu$ de $S(\fa) \otimes \C T$, c'est à dire par $\gamma
\in \fa^*$ et un caractère $\mu$ de $T$. Après avoir détaillé
différentes propriétés de ces modules (poids, dual, etc.)  nous
montrons au \S~\ref{sec3} comment leur structure peut se
ramener à celle des modules de la série principale sur une
algèbre de Hecke graduée associée à un système de racines de
type $A_{i-1} \times A_{n-i-1}$ défini par $\mu$
(cf.~Propositions~\ref{prop11} et~\ref{prop12}) et nous en
déduisons un crit\`ere d'irr\'eductibilit\'e pour $M(\gamma
\otimes \mu)$ (cf.~Th\'eor\`eme~\ref{thm13}).  L'étude de
$\H_B(\tk)$ utilise fortement l'existence d'un théorème de
Poincaré-Birkhoff-Witt (PBW) pour cette algèbre; nous obtenons
ce résultat grâce au théorème de PBW qui existe dans $\che(k)$,
cf.~\cite{EG}.

Pour un système de racines de type $D_n$ une étude de l'algèbre
$\D$ montre qu'elle possède des propriétés (presque) analogues
à celles de $\H_B(\tk)$. Cette étude repose sur l'inclusion
d'un système de racines de type $D_n$ dans un sytème de type
$B_n$ et des résultats obtenus précédemment pour ce cas. C'est
l'objet de la derni\`ere section (cf., en particulier, le
Th\'eor\`eme~\ref{thm48}).

\section{Notations et pr\'eliminaires}
\label{sec1}

Soit $W_B$ un groupe de Weyl de type $B_n$ r\'ealis\'e comme
groupe de r\'eflexions dans un espace vectoriel r\'eel $V$ de
dimension $n$ et $R_B$ le syst{\`e}me de racines correspondant
\`a cette donn\'ee.
On note $V^*$ le dual de $V$, $<.\, ,.>$ le crochet de
dualit{\'e}, et $V^*_\C$, $V_\C$ les complexifi{\'e}s
respectifs. L'action de $W_B$ s'{\'e}tend naturellement en une
action par automorphismes sur les alg\`ebres sym\'etriques
$S(V_\C)$ et $S(V_\C^*)$.  Si $\alpha \in R_B$, on d\'esigne
par $s_\alpha$ la r\'eflexion associ{\'e}e et $\alpha\spcheck
\in V^*$ la coracine.

On fixe une base orthonorm{\'e}e
$\{\epsilon_1,\epsilon_2,..,\epsilon_n\}$ de $V$ form\'ee de
racines courtes.  On choisit $S_B=\{\epsilon_i- \epsilon_{i+1},
1\leq i \leq n-1 \, ; \, \epsilon_n\}$ comme ensemble de
racines simples de sorte que $R_B^+=\{\epsilon_i \pm
\epsilon_j, 1 \leq i<j \leq n \, ; \, \epsilon_k , 1 \leq k
\leq n\}$ est l'ensemble des racines positives d\'etermin\'e
par ce choix.

Pour tous $1 \le i,j \le n$ on pose $\alpha_{i,j}
=\epsilon_i-\epsilon_j$. Soit $R_A$ le sous-ensemble de $R_B$
constitu{\'e} des $\pm \alpha_{i,j}$, $1 \leq i<j \leq n$;
c'est un syst{\`e}me de racines de type $A_{n-1}$. On note
$W_A$ le sous-groupe de $W_B$ associ{\'e}, $S_A=R_A \cap
S_B=\{\alpha_i =\epsilon_i-\epsilon_{i+1}, 1 \leq i \leq n-1\}$
l'ensemble des racines simples et $\ell$ la longueur sur $W_A$
d\'efinie par $S_A$. L'action par permutation de $W_A$ sur la
base $\epsilon_1,..,\epsilon_n$ d{\'e}finit une action sur
$\{1,..,n\}$ en posant $w*i =j$ si $w(\epsilon_i)=\epsilon_j$,
ce qui permet d'identifier $W_A$ avec le groupe sym\'etrique
$\mathfrak{S}_n$.

L'alg{\`e}bre du groupe $W_B$ est d\'esign\'ee par $\C
W_B=\C\langle t_w, w \in W_B\rangle$.  On {\'e}crira
$t_i=t_{s_{\epsilon_i}}$, $1 \leq i \leq n$,
$s_{p,q}=s_{\alpha_{p,q}}$, $1 \leq p<q \leq n$.

Soit $T$ le sous-groupe (normal ab\'elien) de $W_B$ engendr\'e
par les $s_{\epsilon_i}$, $1 \le i \le n$.  On a donc $T \simeq
\{\pm 1\}^{n}$ et $W_B= T \rtimes W_A$.  Si $T\spcheck$
d\'esigne le groupe des caract\`eres de $T$ on identifie $\C
T\spcheck$ au dual $(\C T)^*$, en \'etendant un \'el\'ement de
$T\spcheck$ en une forme lin\'eaire sur $\CT$.  Le groupe $W_A$
op\`ere sur $T\spcheck$ par la formule $\twisted{w}{\mu}(x) =
\mu(w^{-1}xw)$ pour tous $w \in W_A, \mu \in T\spcheck, x \in
T$. On notera $W_A(\mu)$ le stabilisateur de $\mu \in
T\spcheck$ pour cette action.

On rappelle que tout caract\`ere de l'alg\`ebre $S(V_\C)$ est
d\'etermin\'e par un \'el\'ement de $V_\C^*$ et que l'action de
$w \in W_A$ sur $\gamma \in V^*_\C$ est d\'efinie par
$\twisted{w}{\gamma }(\zeta)=\gamma(w^{-1}(\zeta))$ pour tout
$\zeta \in V_\C$.

\subsection{ Alg\`ebre de Hecke gradu\'ee}
\label{sec1.1}

Soit $R \subset V$ un syst\`eme de racines de base $S$, de
groupe de Weyl $W$.  Une multiplicit\'e sur $R$ est une
fonction $W$-invariante $c: R \to \C$.  On peut lui associer
une alg\`ebre de Hecke gradu\'ee $\H_{\mathit{gr}}(c,S)$
d\'efinie de la mani\`ere suivante.  L'alg\`ebre
$\H_{\mathit{gr}}(c,S)$ est engendr\'ee par les \'el\'ements
$\zeta \in V$ et $t_w, w \in W$, soumis aux relations:
\begin{itemize}
\item $\C\langle \zeta \in V \rangle \simeq S(V_\C)$,
  $\C\langle t_w, w \in W \rangle \simeq \C W$;
\item $t_{s_\alpha}\zeta=
  s_\alpha(\zeta)t_{s_\alpha}+<\alpha\spcheck,s_\alpha(\zeta)>
  c_\alpha$ pour tout $\alpha \in S$.
\end{itemize}
Cette alg\`ebre v\'erifie un th\'eor\`eme de
Poincar\'e-Birkhoff-Witt
\begin{equation*}
\tag{$\mathrm{PBW}$0} \H_{\mathit{gr}}(c,S) \simeq
S(V_\C)\otimes\CW 
\end{equation*}
comme espace vectoriel, ou $S(V_\C)$-module (cf.~\cite{Lu}).

Pour tout caract\`ere $\gamma : S(V_\C) \to \C$, d\'efini par
$\gamma \in V_\C^*$, on note $\C v_\gamma$ la repr\'esentation
de dimension~$1$ correspondante et l'on construit un module
induit, appel\'e module de la s\'erie principale,
$$
M(\gamma) =
\mathrm{Ind}_{S(V_\C)}^{\H_{\mathit{gr}}(c,S)}(\C v_\gamma).
$$
On sait qu'alors $M(\gamma) = \bigoplus_{w \in W} \C t_w
\otimes v_\gamma$ et l'on dispose d'un crit\`ere
d'irr\'eductibilit\'e pour ce module (cf.~\cite{KR}).  On se
propose de g\'en\'eraliser cette notion, et le crit\`ere
d'irr\'eductibilit\'e, \`a une nouvelle alg\`ebre associ\'ee
\`a un syst\`eme de racines de type $B$.

\subsection{Une g\'en\'eralisation de l'alg\`ebre de Hecke
  gradu\'ee}
\label{sec1.2}
 
On reprend les notations pr\'ec\'edentes relatives \`a un
syst\`eme de racines de type $B_n$. Soit $\bk= (k,k_c) : R_B
\to \C$ une multiplicit\'e de valeur $k \ne 0$ sur les racines
longues et $k_c$ sur les racines courtes.  G\'en\'eralisons la
notion de multiplicit\'e en d\'efinissant dans ce cadre une
fonction $\tk : R_A \to \CT$ en posant:
$$
\tk_{\alpha_{i,j}}=k(1+t_it_{j})
$$
pour tous $1 \le i \ne j \le n$.  Il r\'esulte de
$wt_it_jw^{-1} = t_{w*i}t_{w*j}$, pour $w \in W_A$, que
$$
\text{$\tk_{w(\alpha)} = w\tk_\alpha w^{-1}$ pour tous $w
  \in W_A$, $\alpha \in R_A$}.
$$
Donc $\tk$ est $W_A$-\'equivariante. Soit alors
$\H_B=\H_B(\tk)$ l'alg\`ebre engendr\'ee par les \'el\'ements
$\zeta \in V$, $t_w \in W_B$ soumis aux relations:

\begin{enumerate}
\item[(a)] $\C\langle \zeta \in V \rangle \simeq S(V_\C)$,
  $\C\langle t_w, w \in W_B \rangle \simeq \C W_B$;
\item[(b)] $t_w \zeta = \zeta t_w$ pour tous $w \in T$, $\zeta
  \in V$;
\item[(c)] $t_{s_\alpha}\zeta=
  s_\alpha(\zeta)t_{s_\alpha}+<\alpha\spcheck,s_\alpha(\zeta)>
  \tk_\alpha$ pour tout $\alpha \in S_A$.
\end{enumerate}
Observons qu'en prenant $\zeta= \epsilon_j$ et $s_\alpha =
s_{\alpha_i}$ dans (c) on obtient:
\begin{equation}
\label{eq1}
t_{s_{\alpha_i}} \epsilon_j = \epsilon_{s_{\alpha_i}*j}
t_{s_{\alpha_i}} - 
<\alpha_i\spcheck,\epsilon_j > \tk_{\alpha_i}.
\end{equation}
Remarquons \'egalement que $\H_B$ est engendr\'ee par les $t_w$
et les $\epsilon_j$ et que, gr\^ace \`a \eqref{eq1}, tout
\'el\'ement de $\H_B$ est combinaison lin\'eaire de mon\^omes
de la forme $\epsilon_1^{j_1} \cdots \epsilon_n^{j_n}t_w$.

\subsection{PBW}
\label{sec1.3}

Soit $\aalg_\mathbb{R}$ un $\mathbb{R}$-espace vectoriel
euclidien de dimension $n$, $R_B \subset \aalg_\mathbb{R}^*$ un
syst\`eme de racines de type $B_n$, $\fa= \C \otimes_\R
\fa_\R$, $\calP=S(\fa^*)$, $\calS=S(\fa)$.  Rappelons qu'\`a la
multiplicit\'e $\bar{k} : R_B \to \C$ on peut associer une
alg{\`e}bre de Cherednik rationnelle not\'ee $\che_B(\bk)$,
$\che(\bk)$ ou $\che$, dont une r\'ealisation dans $\End_\C
\calP$ est fournie, via les op\'erateurs de Dunkl, de la
fa\c{c}on suivante (cf.~\cite{TL,EG}). Si $w \in W_B$ on note
ici $t_w \in \Aut \calP$ l'automorphisme associé et pour tout
$y \in \fa$ on pose
$$
T_y=T_y(\bk) =T_y^B(\bk) = \partial_y + \half\sum_{\alpha
  \in R_B} \bk_\alpha \frac{<{\alpha},{y}>}{\alpha} (1 -
t_{s_\alpha}).
$$
Alors, $\che$ est la sous-alg{\`e}bre de $\End_\C(\calP)$
engendr{\'e}e par les $t_w, w \in W_B$, $x \in \aalg^*$ et les
$T_y$, $y \in \aalg$.  Ces g{\'e}n{\'e}rateurs v\'erifient les
relations suivantes:
\begin{enumerate}   
\item $[T_{y},x] = \dual{y}{x} +\frac{1}{2}\sum_{\alpha \in
    R_B} \bk_{\alpha}<{y},{\alpha}> <{\alpha\spcheck},{x}>
  t_{s_\alpha}$;
\item $t_wxt_{w^{-1}}=w(x)$;
\item $t_wT_{y}t_{w^{-1}}=T_{w(y)}$.
\end{enumerate}   
Il r\'esulte de~\cite[Theorem~1.3]{EG} que l'on a un
th\'eor\`eme de Poincar\'e-Birkhoff-Witt dans $\che$:

\begin{Theoreme}[PBW1]   
\label{PBW1}
Le $\calP$-module $\che$ est isomorphe {\`a}
$\calP\otimes\calS\otimes \C W_B$.
\end{Theoreme}   

Nous allons montrer que $\H_B$ se r\'ealise comme une
sous-alg{\`e}bre de $\che$ \`a l'aide d'op\'erateurs introduits
par S.~Kakei dans~\cite{Ka}.  On fixe une base orthonorm\'ee
$\{e_1,\dots,e_n\}$ de $\fa_\R$, de base duale
$\{z_1,\dots,z_n\}$. Soit $\D_B$ la sous-alg{\`e}bre de $\che$
engendr{\'e}e par les $z_j T_{e_j}$, $j=1,\dots,n$, et les
$t_w$, $w \in W_B$.  On pose:
$$
D_j = D_j^B =z_j T_{e_j} + \sum_{1 \le i < j }
\tk_{\alpha_{i,j}} t_{s_{\alpha_{ij}}}.
$$
(Les $D_j$ sont les op\'erateurs $\DK_{j}^B$ d\'efinis en
\cite[\S~2.3]{Ka}.) Il est clair que
$$
\D_B=\C\langle D_1,\dots,D_n \, ; \, t_w, \, w \in
W_B\rangle.
$$
On v\'erifie facilement que les $D_j$ sont lin\'erairement
ind\'ependants et que l'on a les relations suivantes dans
$\che$:
\begin{enumerate}
\item $D_iD_j=D_jD_i$;
\item $D_j t_i=t_iD_j$;
\item $t_{s_{\alpha_i}}D_j = D_{s_{\alpha_i}*j}t_{s_{\alpha_i}}
  - <\alpha_i\spcheck,\epsilon_j > \tk_{\alpha_i}$.
\end{enumerate}
On pose $V'=\bigoplus_{j=1}^n \R D_j$; comme les $D_j$
commutent, la $\C$-alg\`ebre qu'ils engendrent est un quotient
de $S(V'_\C)$. Elle est en fait \'egale \`a $S(V'_\C)$, comme
le montre le r\'esultat suivant:

\begin{Lemme}[PBW2]
\label{PBW2}
Le $\C$-espace vectoriel $\D_B$ est isomorphe \`a $S(V'_\C)
\otimes \C W_B$.
\end{Lemme}

\begin{proof}
  Si $j=(j_1,\dots,j_n) \in \N^n$ on posera $|j| =j_1 + \cdots
  + j_n$, $z^j =z_1^{j_1}\cdots z_n^{j_n}$, $T^j =
  T_{e_1}^{j_1}\cdots T_{e_n}^{j_n}$ et $D_j=D_1^{j_1}\cdots
  D_n^{j_n}$. Il r\'esulte des relations ci-dessus que tout
  \'el\'ement de $\D_B$ peut s'\'ecrire
  $$
  \textstyle{P=\sum_{d=0}^m \sum_{|j|=d,\, w \in W_B}
    \lambda_{j,w}D^j t_w }
  $$
  pour des $\lambda_{j,w} \in \C$; il s'agit de voir que
  cette \'ecriture est unique. Supposons $P=0$.  Il existe
  alors un couple $(j,w)$ tel que $\lambda_{j,w} \neq 0$ et
  $|j|=m$ maximal pour cette propri\'et\'e.  La d\'efinition
  des $D_p$ et les relations de commutation dans $\che$
  entra\^{\i}nent que si $|p|=d$, $D^p= z^pT^p + Q$ avec $Q \in
  \calP\calS_{d-1}\C W_B$, o\`u l'on a pos\'e
  $\calS_{d-1}=\bigoplus_{|i| \le d-1} \CT^i$. On obtient ainsi
  $P=\sum_{|j|=m, \, w \in W_B} \lambda_{j,w} z^jT^j t_w + R$
  avec $R \in \calP \calS_{m-1}\C W_B$. Le
  Th\'eor\`eme~\ref{PBW1} implique alors que $\lambda_{j,w}=0$
  pour tout $|j|=m$, ce qui contredit l'hypoth\`ese.
\end{proof}

\begin{Theoreme}[PBW3]
\label{PBW3}
L'application $\Phi : \H_B \to \D_B$ d\'efinie par
$$
\Phi(\epsilon_j)= D_j, \quad \Phi(t_w) = t_w, \ \,
\text{pour $j=1,\dots,n$ et $w \in W_B$,}
$$
est un isomorphisme d'alg\`ebres. En particulier, $\H_B$ est
isomorphe \`a $S(V_\C) \otimes \C W_B$ comme $S(V_\C)$-module.
\end{Theoreme}

\begin{proof}
  La d\'efinition par g\'en\'erateurs et relations de $\H_B= \C
  \langle \epsilon_1,\dots,\epsilon_n \, ; \, t_w, \, w \in W_B
  \rangle$ et les relations pr\'ec\'edentes dans $\D_B$
  montrent que $\Phi$ d\'efinit un morphisme surjectif
  d'alg\`ebres. Le Lemme~\ref{PBW2} assure que $\Phi$ est
  injectif. La deuxi\`eme assertion est alors \'evidente.
\end{proof}

\subsection{Relations de commutation dans $\H_B$} 
\label{sec1.4}

Pr\'ecisons les relations de commutation dans $\H_B$.
Rappelons que l'ensemble d'inversion de $w^{-1} \in W_A$ est
$R(w^{-1}) = \{ \alpha \in R^+_A : w^{-1}(\alpha) \in R^-_A\}$
et que $\ell(s_\alpha w) < \ell(w)$ pour tout $\alpha \in
R(w^{-1})$.

\begin{Lemme}
\label{lem1}
Pour tous $\zeta \in V$ et $w \in W_A$, on a:
$$
t_w \zeta = w(\zeta)t_w + \sum_{\alpha \in R(w^{-1})}
<\alpha\spcheck,w(\zeta)>\tk_{\alpha}t_{s_{\alpha} w} =
w(\zeta)t_w + \sum_{\alpha \in R(w^{-1})}
<\alpha\spcheck,w(\zeta)> t_{s_{\alpha} w}
\tk_{w^{-1}(\alpha)}.
$$
En particulier, pour $1 \leq p<q \leq n$ :
\begin{equation*}
\begin{split}
  t_{s_{p,q}}\zeta=s_{p,q}(\zeta)t_{s_{p,q}}+ & \sum_{p+1 \leq
    j
    <q}<\alpha_{j,q}\spcheck,s_{p,q}(\zeta)>t_{s_{j,q}s_{p,q}}
  \tk_{\alpha_{j,p}}+ \sum_{p+1 \leq
    j<q}<\alpha_{p,j}\spcheck,s_{p,q}(\zeta)>t_{s_{p,j}s_{p,q}}
  \tk_{\alpha_{j,q}} \\ &
  +<\alpha_{p,q}\spcheck,s_{p,q}(\zeta)>\tk_{\alpha_{p,q}}.
\end{split}
\end{equation*}
\end{Lemme}  

\begin{dem}  
  La multiplicit\'e $\tk$ \'etant $W_A$-\'equivariante et
  commutant aux \'el\'ements de $V$, la premi\`ere formule se
  d\'emontre comme celle de \cite[Proposition~1.1(1)]{Op2} dans
  une alg\`ebre de Hecke gradu\'ee. La seconde r\'esulte de
  $R(s_{p,q}^{-1})=R(s_{p,q})=\{\alpha_{j,q},\, \alpha_{p,i} :
  p \leq j <q, \, p+1 \leq i <q \}$.
\end{dem}

D\'efinissons des << op{\'e}rateurs de diff{\'e}rences
divis{\'e}es >> $\Delta_j \in \End S(V_\C)$ par
$$
\Delta_j(p)=k (p-s_{\alpha_j}(p))/{\alpha_j}
$$
pour $p \in S(V_\C)$ et $j \in \{1,\dots,n-1\}$. Rappelons
que dans $\H_{gr}(k,S_A)$ on a $\Delta_j(p) =
t_{s_{\alpha_j}}p-s_{\alpha_j}(p)t_{s_{\alpha_j}}$.  On
g\'en\'eralise donc ces op\'erateurs en d{\'e}finissant, pour
tout $j \in \{1,\dots,n-1\}$, $\tdel_j: S(V_\C) \rightarrow
\H_B$ par
$$
\tdel_j(p)=t_{s_{\alpha_j}}p -
s_{\alpha_j}(p)t_{s_{\alpha_j}},
$$
o\`u ici le calcul s'effectue dans $\H_B$.

\begin{Lemme} 
\label{lem2}
Soit $p \in S(V_\C)$, alors pour tout $j \in \{1,\dots,n-1\}$,
on a:
$$
\tdel_j(p)=\Delta_j(p) (1+t_jt_{j+1}) = \tk_{\alpha_j}
(p-s_{\alpha_j}(p))/{\alpha_j}.
$$
\end{Lemme} 
 
\begin{dem} 
  D\'emontrons le résultat par r\'ecurrence sur le degr\'e de
  $p$.  Pour $p$ de degr\'e $1$, il provient directement des
  relations (c), \S~\ref{sec1.2}, dans $\H_B$.  Supposons qu'il
  soit vrai pour $p_1,p_2 \in S(V_\C)$ et montrons le pour le
  produit $p_1p_2$.  On a:
\begin{equation*}
\begin{split}
  \tdel_j(p_1p_2) & =
  t_{s_{\alpha_j}}p_1p_2-s_{\alpha_j}(p_1)s_{\alpha_j}(p_2)
  t_{s_{\alpha_j}}
  \\
  & =
  (t_{s_{\alpha_j}}p_1-s_{\alpha_j}(p_1)t_{s_{\alpha_j}})p_2 +
  s_{\alpha_j}(p_1)\tdel_j(p_2) \\
  & = \tdel_j(p_1)p_2 + s_{\alpha_j}(p_1)\tdel_j(p_2) \\ &
  =\Delta_j(p_1)p_2(1+t_jt_{j+1})+s_{\alpha_j}(p_1)
  \Delta_j(p_2)(1+t_jt_{j+1}),
\end{split}
\end{equation*}
car les \'el\'ements de $T$ et $S(V_\C)$ commutent.  Il suffit
alors pour conclure de remarquer que
$\Delta_j(p_1)p_2+s_{\alpha_j}(p_1)\Delta_j(p_2)=
\Delta_j(p_1p_2)$.
\end{dem} 
 
Rappelons que $k\ne 0$. De (PBW3) il r\'esulte donc que
$$
\tdel_j(p) = 0 \iff \Delta_j(p) = 0 \iff s_{\alpha_j}(p) =p.
$$
Ainsi l'intersection des noyaux (dans $S(V_\C)$) des
$\tdel_j$, $j=1,\dots,n$, est l'alg\`ebre $S(V_\C)^{W_A}$ des
\'el\'ements $W_A$-invariants de $S(V_\C)$.

\subsection{Centre de $\H_B$}
\label{sec1.5}
   
Rappelons que tout \'el\'ement $t_x$, $x \in T$, est de la
forme $t_{\underline{i}}=t_{i_1}\cdots t_{i_j}$ pour un
multi-indice $\underline{i}=\{1 \le i_1 < \cdots <i_j \le n\}$
de longueur $j \in \{0,\dots,n\}$. Remarquons que le groupe
$W_A$ op\`ere dans $\CT$ par conjugaison et que
$t_{\underline{i}}$ est conjugu\'e (par cette action) \`a
$t_{\underline{i}'}$ si et seulement si $|\underline{i}|=
|\underline{i}'|$.  Pour chaque $j \in \{1,..,n\}$, on pose:
$$
\vt_j = \sum_{|\underline{i}|=j} t_{\underline{i}}.
$$
Il est alors facile de voir que l'alg\`ebre d'invariants
$(\CT)^{W_A}$ est \'egale \`a $\C[\vt_1, \dots,
\vt_n]=\bigoplus_{j=0}^n \C \vt_j$.

\begin{lem}
\label{lem3}
Le centralisateur $\{c\in \H_B : [c,V] = 0 \}$ de $V$ dans
$\H_B$ est la sous-alg\`ebre (commutative) $S(V_\C) \otimes
\CT$.
\end{lem}

\begin{dem}
  Soit $0 \ne c \in \H_B$; par (PBW3) on peut {\'e}crire de
  mani{\`e}re unique $c$ sous la forme $\sum_{w \in W_A, \, x
    \in T} p_{w,x}t_wt_x$ pour des $p_{w,x} \in S(V_\C)$ non
  tous nuls.  Posons $\calC=\{(w,x) \in W_A \times T: \;
  p_{w,x} \ne 0\}$ et notons $m = \max\{\ell(w) :\exists x \in
  T, \; (w,x) \in \calC\}$. Si $\xi \in V$ il vient, en
  utilisant le Lemme~\ref{lem1},
\begin{equation*}
\begin{split}
  [c,\xi] & = \textstyle{\sum_{(w,x) \in \calC}} \, p_{w}
  (t_w\xi - \xi t_w)t_x \\
  & = \textstyle{\sum_{(w,x) \in \calC, \, \ell(w)=m}} \, p_{w}
  (w(\xi) - \xi) t_wt_x + \textstyle{\sum_{(w',x') \in W_A
      \times T, \, \ell(w') < m}} \, q_{w',x'}t_{w'}t_{x'}
\end{split}
\end{equation*}
avec $q_{w',x'} \in S(\VC)$. Par cons\'equent (PBW3) assure que
si $c$ est dans le centralisateur de $V$ et $(w,x) \in \calC$
v\'erifie $\ell(w)=m$, alors $w(\xi) = \xi$ pour tout $\xi$,
i.e.~$w=\mathrm{id}$. Donc $m=0$ et $c \in S(V_\C) \otimes
\CT$.
\end{dem}

Rappelons que le centre d'une alg\`ebre de Hecke gradu\'ee
$\H_{\mathit{gr}}(c,S)$ comme en~\ref{sec1.1} est $S(\VC)^{W}$.
Pour $\H_B$ le centre, $Z(\H_B)$, est donn\'e par le
th\'eor\`eme suivant.

\begin{Theoreme} 
\label{thm4}
On a $Z(\H_B)=S(V_\C)^{W_A}\otimes(\CT)^{W_A}=\bigoplus_{j=0}^n
S(V_\C)^{W_A}\vt_j.$
\end{Theoreme}   
 
\begin{dem} 
  Observons que $c \in Z(\H_B)$ si et seulement si
  $[c,V]=[c,T]=[t_{s_{\alpha_j}},c] = 0$ pour $j=1,\dots,n-1$.
  Par le Lemme~\ref{lem3}, ceci \'equivaut \`a $c \in S(\VC)
  \otimes \CT$ et $[t_{s_{\alpha_j}},c] = 0$ pour
  $j=1,\dots,n-1$. \'Ecrivons $c=\sum_{x \in T} p_x t_x$ avec
  $p_x \in S(\VC)$. Avec les notations de~\ref{sec1.4} on
  trouve que
\begin{equation*}
\begin{split}
  [t_{s_{\alpha_j}},c] &= \sum_x t_{s_{\alpha_j}}p_x t_x
  -\sum_x
  p_x t_{s_{\alpha_j}} t_{s_{\alpha_j}xs_{\alpha_j}}  \\
  & = \sum_x \bigl(s_{\alpha_j}(p_x) -
  p_{s_{\alpha_j}xs_{\alpha_j}}\bigr) t_{s_{\alpha_j}} t_{x} +
  \sum_x \tdel_j(p_x) t_x.
\end{split}
\end{equation*}
On d\'eduit donc de (PBW3) que $[t_{s_{\alpha_j}},c]=0$ si, et
seulement si, $\tdel_j(p_x)=0$ et
$s_{\alpha_j}(p_x)=p_{s_{\alpha_j}xs_{\alpha_j}}$ pour tous $x
\in T$, $1 \le j \le n-1$. Par la remarque suivant le
Lemme~\ref{lem2} ces conditions \'equivalent \`a $p_x \in
S(\VC)^{W_A}$ et $p_x= p_{wxw^{-1}}$ pour tout $x \in T$, ce
qui compte tenu de la d\'efinition des $\vt_j$ s'\'ecrit aussi
$c= \sum_{j=0}^n p_j \vt_j$ avec $p_j \in S(V_\C)^{W_A}$.
\end{dem}

\subsection{Automorphismes et anti-automorphismes de $\H_B$}
\label{sec1.6}
Soit $R$ une $\C$-algèbre. On note $\RMod{R}$,
resp.~$\ModR{R}$, la catégorie des $R$-modules à gauche,
resp.~droite. La sous-catégorie de $\RMod{R}$ formée des
modules de dimension finie est désignée par $\Rmod{R}$. Le
treillis, ordonné par inclusion, des sous-modules de $M \in
\RMod{R}$ sera noté $\calL_R(M)$ et la longueur de $M$ par
$\lg_R(M)$.  Soit $S$ une autre $\C$-algèbre et $N \in
\RMod{S}$; on écrira $\calL_R(M) \cong \calL_S(N)$ pour
signifier qu'il existe un isomorphisme de $\C$-espaces
vectoriels $f: M \to N$ tel que $X \mapsto f(X)$ induise un
isomorphisme entre les treillis ordonnés $\calL_R(M)$ et
$\calL_S(N)$.  Tout $M \in \Rmod{R}$ est de longueur finie et
on peut lui associer la suite $\calS_R(M)=\{[S_i(M)]\}_{1\le i
  \le t}$ des classes d'isomorphisme des sous-quotients simples
de $M$, comptés sans multiplicité.  Il est clair que si
$\calL_R(M) \cong \calL_S(N)$, une suite de composition de $M$
est transformée par $f$ en une suite de composition de $N$ et
les $[f(S_i(M))]$ sont donc les éléments de la suite
$\calS_S(N)$.

Rappelons que si $M \in \RMod{R}$ et $\vkappa \in \Aut_\C(R)$,
on peut définir un module $\twisted{\vkappa}{M} \in \RMod{R}$
de la façon suivante: $\twisted{\vkappa}{M} = M$ comme
$\C$-espace vectoriel, muni de l'action $a\centerdot u
=\vkappa(a).u$ pour $a \in R, u \in M$.  Notons que cette
torsion par $\vkappa$ laisse stable $\Rmod{R}$ et conserve la
longueur. Rappelons aussi qu'un anti-automorphisme
$\C$-linéaire $\iota$ de $R$ fournit un isomorphisme entre
$\RMod{R}$ et $\ModR{R}$: on fait de $M \in \RMod{R}$ un
$R$-module à droite en posant $M^{\iota} = M$ comme $\C$-espace
vectoriel, muni de la structure de $R$-module à droite définie
par $u.a = \iota(a).u$. L'existence de $\iota$ permet de munir
le dual $M^* = \Hom_\C(M,\C)$ d'une structure de $R$-module à
gauche par la formule
$$
< a.f, u > \, = \, < f, \iota(a).u >,
$$
pour tous $a \in R, f \in M^*, u \in M$. Si de plus $\iota$
est involutif et $M \in \Rmod{R}$ l'application bijective
canonique $c_M : M \to M^{**}$ est alors un isomorphisme de
$R$-modules et l'on a $\lg_R(M) = \lg_R(M^*)$.

Donnons des exemples dans l'algèbre $\hb$ pour lesquels on peut
effectuer les constructions précédentes. Observons qu'un
élément de $\GL(\hb)$ définit un (anti)-automorphisme si, et
seulement si, il préserve les relations (a), (b) et (c)
introduites en~\ref{sec1.2}.

\noindent{(1)} Chaque   $w \in W_B$ définit un automorphisme
intérieur de $\hb$ donné par $\Int(w) : a \mapsto t_w a
t_{w^{-1}}$; le module tordu par l'action de $\Int(w)$ sera
noté $\twisted{w}{M}$.

\noindent{(2)} Le déterminant dans $\GL(V)$
définit le caractère $\sgn : W_B \to \{\pm 1\}$ par $\sgn(w) =
\det(w)$. En particulier, $\sgn$ fournit un élément de
$T\spcheck$ et pour $\mu \in T\spcheck$ on pose $-\mu =\sgn
\otimes \mu \in T\spcheck$; ce caractère est déterminé par
$(-\mu)(t_j) = -\mu(t_j)$ pour tout $j=1,\dots,n$.
L'application $\delta : x \mapsto \sgn(x) t_x$, $x \in T$,
induit un automorphisme de $\CT$ tel que
$\delta(t_it_j)=t_it_j$ pour tous $i,j$. On obtient un
automorphisme de $\hb$ en étendant $\delta$ comme suit:
$$
\delta(\xi) = \xi, \quad \delta(t_x) = \sgn(x)t_x, \quad
\delta(t_w) = t_w,
$$
pour tous $\xi \in V, x \in T, w \in W_A$.

\noindent{(3)} On vérifie facilement  que les formules
ci-dessous donnent un anti-automorphisme involutif $\iota$
de~$\hb$:
$$
\iota(\xi) = -\xi, \quad \iota(t_w) = \sgn(w)t_{w^{-1}},
$$
pour tous $\xi \in V, w \in W_B$.

Il résulte en particulier de (3) que $\RMod{\hb} \simeq
\ModR{\hb}$ et $\Rmod{\hb} \simeq \modR{\hb}$.

\section{Modules de la s{\'e}rie principale}
\label{sec2}

On va d\'efinir pour $\H_B$, de m\^eme que pour une alg\`ebre
de Hecke gradu\'ee, les modules de la s\'erie principale pour
lesquels on démontrera dans la section suivante un crit\`ere
d'irr\'eductibilit\'e.

\subsection{Poids des $\H_B$-modules de dimension finie}   
\label{sec2.1}
 
Un morphisme d'alg\`ebres $\tgamma : S(V_\C)\otimes \CT
\rightarrow \C$ est uniquement déterminé par une forme
lin\'eaire $\gamma \in V_\C^*$ et un caract\`ere $\mu \in
T\spcheck$; nous le noterons $\gamma\otimes\mu$.  On a donc
$(\gamma\otimes\mu)(v \oplus t)=\gamma(v)\mu(t)$, pour $v \in
V_\C$, $t \in \CT$.  On désigne par $\calA =\{\gamma\otimes \mu
: \gamma \in V_\C^*, \, \mu \in T\spcheck\}$ l'ensemble des
caractères de $S(V_\C)\otimes\CT$.

Soient $M$ un $\H_B$-module de dimension finie et $\tnu \in
\calA$. Les sous-espaces poids et sous-espaces poids
g\'en\'eralis\'es de $M$ associ\'es \`a $\tnu$ sont
respectivement:
$$
M_{\tnu}=\{ m \in M : \forall a \in V_\C\oplus\CT, \,
a.m=\tnu(a)m \}
$$
et
$$
M_{\tnu}^{\mathit{gen}}=\{ m \in M : \forall a \in
V_\C\oplus\CT, \, \exists k \in \N, (a-\tnu(a))^k.m=0 \}.
$$
On dira que $\tnu$ est un {poids} de $M$ si $M_{\tnu}^{gen}
\ne 0$ (ce qui équivaut à $M_\tnu \ne 0$).  Comme
$S(V_\C)\otimes\CT$ est une sous-alg\`ebre commutative de
$\H_B$, on a:
$$
M=\bigoplus_{\tnu \in \calA}M_{\tnu}^{\mathit{gen}}.
$$

\subsection{Modules de la série principale}
\label{sec2.2}

Pour tout $\tgamma =\gamma \otimes \mu \in \calA$ on notera
$\twisted{w}{\tgamma} = \twisted{w}{\gamma} \otimes
\twisted{w}{\mu}$ l'action diagonale de $w \in W_A$ et on
désigne par $\C v_\tgamma$ le $S(V_\C)\otimes\CT$-module de
dimension $1$ d\'efini par: $a.v_\tgamma = \tgamma(a)v_\tgamma$
pour tout $a \in S(V_\C)\otimes \CT$.  Le $\H_B$-{module de la
  s\'erie principale} $M(\tgamma)$ associé à $\tgamma$ est le
module induit de $\C v_\tgamma$ \`a $\H_B$:
$$
M(\tgamma)=\H_B\otimes_{S(V_\C^*)\otimes\CT}\C v_\tgamma.
$$
Ce module admet clairement pour base $\{t_w \otimes
v_{\tgamma} : w \in W_A\}$, que l'on peut ordonner de fa\c{c}on
compatible \`a la longueur dans $W_A$.  D'apr\`es les formules
de commutation dans $\H_B$, pour tous $\zeta \in V$, $x \in T$
et $w \in W_A$, il vient:
\begin{equation*}
\label{eqbase}
\zeta t_x.t_w \otimes v_{\tgamma} = \twisted{w}{\tgamma}(\zeta 
x) t_w \otimes v_{\tgamma}-\sum_{\alpha \in R(w^{-1})}
\twisted{w}{\mu}(\tk_{\alpha} 
x)<\alpha\spcheck,\zeta>t_{s_{\alpha}w} \otimes v_{\tgamma}.
\end{equation*}
Comme $\ell(s_{\alpha}w)<\ell(w)$ pour tout $\alpha \in
R(w^{-1})$ les \'el\'ements de $V_\C \bigoplus \CT$ sont
simultan\'ement trigonalis\'es dans la base $\{t_{w} \otimes
v_{\tgamma} : w \in W_A\}$, et les poids de $M(\tgamma)$ sont
les $\{\twisted{w}{\tgamma} : w \in W_A\}$.

La propri\'et\'e suivante montre que l'\'etude des
$\H_B$-modules simples se ram\`ene \`a l'\'etude des modules de
la s\'erie principale sur $\H_B$.

\begin{Prop}
\label{prop5}
Soit $M$ un $\H_B$-module irr\'eductible et $\tgamma$ un poids
de $M$. Alors $M$ est un quotient de $M(\tgamma)$.  En
particulier, un $\H_B$-module simple est de dimension au plus
$|W_A|$.
\end{Prop}

\begin{dem}
  Soit $m_{\tgamma} \in M_{\tgamma} \smallsetminus \{0\}$;
  alors $\C m_{\tgamma}$ est un $S(V_\C)\otimes \CT$-module
  irr\'eductible isomorphe \`a $\C v_\tgamma$.  Comme
  l'induction est le foncteur adjoint de la restriction, il
  existe un unique morphisme de $\H_B$-modules de $M(\tgamma)$
  dans $M$ envoyant $1 \otimes v_\tgamma$ sur $m_\tgamma$.
  Puisque $M$ est irr\'eductible ce morphisme est surjectif et
  $M$ est un quotient de $M(\tgamma)$.
\end{dem}

Rappelons que l'on peut munir le dual de $M \in \Rmod{\hb}$
d'une structure de $\hb$-module grâce à l'anti-automorphisme
$\iota$, et que l'on peut << tordre >> $M$ par $\delta$,
cf.~\ref{sec1.6}. Nous désignerons par $w_0$ l'élément de plus
grande longueur de $W_A$. Pour tout $\tgamma=\gamma \otimes \mu
\in \calA$ on pose
$$
\tgamma^* = (-\twisted{w_0}{\gamma}) \otimes
(-\twisted{w_0}{\mu}).
$$

\begin{Prop}
\label{prop6}
Soit $\tgamma=\gamma \otimes \mu \in \calA$. Alors:
$$
\twisted{\delta}{M(\tgamma)} \simeq M(\gamma \otimes
(-\mu)), \qquad M(\tgamma)^* \simeq M(\tgamma^*).
$$
Tout $\H_B$-module irr\'eductible est un sous-module d'un
module de la série principale.
\end{Prop}

\begin{dem}
  Il est clair que $\twisted{\delta}{M(\tgamma)}$ est engendré
  par $1 \otimes v_{\tgamma}$. De plus, $\delta(\xi).1 \otimes
  v_{\tgamma} = \xi.1 \otimes v_{\tgamma} = \gamma(\xi).1
  \otimes v_{\tgamma}$ et $\delta(t_x).1 \otimes v_{\tgamma}=
  \sgn(x)t_x. 1 \otimes v_{\tgamma}= \sgn(x)\mu(x).1 \otimes
  v_{\tgamma} = (-\mu)(x).1 \otimes v_{\tgamma}$, pour $\xi \in
  V, x \in T$. La propriété universelle de $M(\gamma \otimes
  (-\mu))$ fournit alors une surjection $M(\gamma \otimes
  (-\mu)) \to \twisted{\delta}{M(\tgamma)}$, qui est un
  isomorphisme (puisque ces modules sont de même dimension).
  
  Pour tout $w \in W_A$, notons $f_{w,\tgamma} \in
  M(\tgamma)^*$ la forme linéaire définie par
  $$
  <f_{w,\tgamma}, t_u \otimes v_{\tgamma} > =
\begin{cases}
  \sgn(w) & \ \text{si $u=w$}; \\
  0 & \ \text{sinon.}
\end{cases}
$$
On a donc $M(\tgamma)^* = \bigoplus_{w \in W_A} \C
f_{w,\tgamma}$ et un calcul facile montre que
$t_g.f_{w,\tgamma}= f_{gw,\tgamma}$ pour tout $g \in W_A$. Il
en découle que $M(\tgamma)^* = \hb. f_{w_0,\tgamma}$. Si l'on
montre que $f_{w_0,\tgamma}$ est de poids $\tgamma^*$, on
obtiendra une surjection $M(\tgamma^*) \to M(\tgamma)^*$ qui
donnera l'isomorphisme voulu.

Soient $\xi \in V, w \in W_A$; observons que le
Lemme~\ref{lem1} implique $-\xi.t_{w}\otimes v_\tgamma =
-\gamma(w^{-1}(\xi)).t_{w} \otimes v_\tgamma + \sum_{u \in W_A,
  \, \ell(u) < \ell(w)} \lambda_u.t_u \otimes v_\tgamma$ avec
$\lambda_u \in \C$. On a donc $< \xi.f_{w_0,\tgamma},
t_{w}\otimes v_\tgamma > = 0$ sauf si $w=w_0$, auquel cas on
trouve $-\gamma(w_0(\xi))\sgn(w_0) =
-\twisted{w_0}{\gamma}(\xi)\sgn(w_0)$. Par conséquent
$\xi.f_{w_0,\tgamma}= -\twisted{w_0}{\gamma}(\xi)
f_{w_0,\tgamma}$.

Soit maintenant $x \in T$. On a $<t_x.f_{w_0,\tgamma},
t_{w}\otimes v_\tgamma > = < f_{w_0,\tgamma}, \sgn(x)t_x.
t_{w}\otimes v_\tgamma >$. Mais de $\sgn(x)t_x. t_{w}\otimes
v_\tgamma = \sgn(x) t_w \otimes (t_{w^{-1}xw}.v_\tgamma) =
\sgn(x)\mu(w^{-1}xw).t_{w}\otimes v_\tgamma =
-\twisted{w}{\mu}(x).t_{w}\otimes v_\tgamma$ on tire que
$<t_x.f_{w_0,\tgamma}, t_{w}\otimes v_\tgamma >$ est nul sauf
pour $w=w_0$, auquel cas on obtient $ -\twisted{w_0}{\mu}(x)
\sgn(w_0)$. Donc $t_x.f_{w_0,\tgamma} = -\twisted{w_0}{\mu}(x)
f_{w_0,\tgamma}$.

La dernière assertion résulte du résultat précédent, de la
Proposition~\ref{prop5}, et des propriétés générales sur la
dualité rappelées en~\ref{sec1.6}.
\end{dem}

\begin{prop}
\label{prop61}
Soient $\tgamma=\gamma \otimes \mu \in \calA$ et $w \in W_A$.
Il existe un isomorphisme de $\H_B$-modules
$$
M(\tgamma) \longisomto \twisted{w}{M(\tgamma)}, \quad t_g
\otimes v_\tgamma \mapsto t_{wg} \otimes v_\tgamma.
$$
\end{prop}

\begin{dem}
  Consid\'erons l'\'el\'ement $t_w \otimes v_\tgamma$ de
  $\twisted{w}{M(\tgamma)}$. Puisque $t_{w^{-1}u}\centerdot t_w
  \otimes v_\tgamma = t_u\otimes v_\tgamma$ pour tout $u \in
  W_A$, cet l'\'el\'ement engendre le $\H_B$-module
  $\twisted{w}{M(\tgamma)}$. Rappelons (cf.~Lemme~\ref{lem1})
  que pour tout $\eta \in V_\C$ on a
\begin{align*}
  \Int(w)(\eta) &= w(\eta) + \sum_{\alpha \in R(w^{-1})} <
  \alpha\spcheck,w(\eta) > \tk_\alpha t_{s_\alpha},  \\
  \intertext{ou encore} \eta t_w &= t_w w^{-1}(\eta) -
  \sum_{\alpha \in R(w^{-1})} < \alpha\spcheck,\eta> \tk_\alpha
  t_{s_\alpha w}.
\end{align*}
Si $\xi \in V_\C$ il vient donc, dans
$\twisted{w}{M(\tgamma)}$,
\begin{align*}
  \xi\centerdot t_w \otimes v_\tgamma &= w(\xi)t_w \otimes
  v_\tgamma + \sum_{\alpha \in R(w^{-1})} <
  \alpha\spcheck,w(\xi)
  > \tk_\alpha t_{s_\alpha w}  \otimes v_\tgamma \\
  &= t_w \xi\otimes v_\tgamma - \sum_{\alpha \in R(w^{-1})} <
  \alpha\spcheck,w(\xi) > \tk_\alpha t_{s_\alpha w} \otimes
  v_\tgamma + \sum_{\alpha \in R(w^{-1})} <
  \alpha\spcheck,w(\xi)
  > \tk_\alpha t_{s_\alpha w} \otimes v_\tgamma \\
  &= \gamma(\xi) t_w \otimes v_\tgamma.
\end{align*}
D'autre part si $x \in T$ on a: $t_x \centerdot t_w \otimes
v_\tgamma = t_w t_x \otimes v_\tgamma = \mu(x) t_w \otimes
v_\tgamma$.  Ainsi $t_w \otimes v_\tgamma\in
\twisted{w}{M(\tgamma)}$ est de poids $\tgamma$. La propriété
universelle de $M(\tgamma)$ fournit un morphisme surjectif de
$M(\tgamma)$ sur le $\H_B$-module $\H_B\centerdot t_w \otimes
v_\tgamma = \twisted{w}{M(\tgamma)}$, d\'etermin\'e par $1
\otimes v_\tgamma \mapsto t_w \otimes v_\tgamma$; par
cons\'equent $t_g \otimes v_\tgamma \mapsto t_g\centerdot t_w
\otimes v_\tgamma = t_{wg} \otimes v_\tgamma$ pour tout $g \in
W_A$. Comme $\dim M(\tgamma) = \dim \twisted{w}{M(\tgamma)}$ ce
morphisme est bijectif.
\end{dem}

Soient $\H_{\mathit{gr}}(c)=\H_{\mathit{gr}}(c,S)$ une
alg\`ebre de Hecke gradu\'ee comme au \S~\ref{sec1.1} et
$M(\lambda)$ le $\H_{\mathit{gr}}(c)$-module de la s\'erie
principale associ\'e \`a $\lambda \in V_\C^*$.  Pour tout $w
\in W$ on peut tordre $M(\lambda)$ par $\Int(w) \in
\Aut_\C(\H_{\mathit{gr}}(c))$ et obtenir le module
$\twisted{w}{M(\lambda})$. Avec ces notations, une preuve
identique \`a celle de la proposition pr\'ec\'edente donne:

\begin{prop}
\label{prop62}
Soient $\lambda \in V_\C^*$ et $w \in W$. Il existe un
isomorphisme de $\H_{\mathit{gr}}(c)$-modules
$$
M(\lambda) \longisomto \twisted{w}{M(\lambda)}, \quad t_g
\otimes v_\lambda \mapsto t_{wg} \otimes v_\lambda.
$$
\end{prop}

\smallskip

\begin{Remarque} 
\label{rem63}
Rappelons \cite[Proposition~2.8]{KR} que si le
$\H_{\mathit{gr}}(c)$-module $M(\lambda)$ est simple, il en est
de m\^eme de $M(\twisted{w}{\lambda})$ pour tout $w \in W$. Si
$\lambda, \gamma \in V_\C^*$ sont quelconques, l'existence d'un
isomorphisme $M(\lambda) \simeq M(\gamma)$ force $\gamma \in
W.\lambda = \{\twisted{w}{\lambda} : w \in W\}$ (car $\gamma$
est alors un poids de $M(\lambda)$). Par cons\'equent sous
l'hypoth\`ese $M(\lambda)$ simple, la condition $M(\lambda)
\simeq M(\gamma)$ est \'equivalente \`a $\gamma \in W.\lambda$.
\end{Remarque}

\section{Crit\`ere d'irr\'eductibilit\'e}
\label{sec3}

\subsection{Notations}
\label{sec3.1}

Pour tout $\tnu \in \calA$ on d\'esignera par $\psi : \H_B
\rightarrow \End M(\tnu)$ la repr\'esentation $\H_B$ dans
$M(\tnu)$. S'il n'y a pas risque d'ambigu\"{\i}t\'e, l'action
de $a \in \hb$ sur $v \in M(\tnu)$ sera simplement not\'ee $a.v
= \psi(a)(v)$. Soit $\varpi \in T\spcheck$; rappelons que
$W_A(\varpi)$ est le stabilisateur de $\varpi$ dans $W_A$. On
notera $\H_B(\varpi)$ la sous-alg\`ebre de $\H_B$ engendr\'ee
par $V$, $W_A(\varpi)$ et $T$.

On fixe un caract\`ere $\tgamma=\gamma\otimes\mu \in \calA$ et
un syst\`eme de repr\'esentants $\{w_1=\id,w_2,..,w_s\}$ de
$W_A/W_A(\mu)$.

\begin{Remarque}
\label{rem7}
Les caract\`eres $\twisted{w_j}{\mu}$, $j=1,\dots,s$ sont deux
\`a deux distincts. Par un argument classique il en d\'ecoule
que les formes lin\'eaires $\twisted{w_j}{\mu} \in (\CT)^*$
sont lin\'eairement ind\'ependantes. On peut ainsi trouver des
\'el\'ements $y_i \in T$, $i=1,\dots,s$, tels que la matrice
$[\twisted{w_j}{\mu}(t_{y_i})]_{1 \le i,j \le s}$ est
inversible.
\end{Remarque}

Pour tout $j \in \{1,..,s\}$, on d\'efinit un sous-espace
vectoriel de $M(\tgamma)$ par:
$$
E_j(\tgamma)=\bigoplus_{w \in W_A(\mu)} \C t_{w_j}t_w\otimes
v_\tgamma.
$$
Observons que
$$
E_1(\tgamma)= \bigoplus_{w \in W_A(\mu)} \C t_w\otimes
v_\tgamma, \qquad E_j(\tgamma)=t_{w_j}.E_1(\tgamma).
$$
On a $M(\tgamma)=\bigoplus_{j=1}^s E_j(\tgamma)$ et cette
d\'ecomposition est la d\'ecomposition isotypique du $T$-module
$M(\tgamma)$: le groupe $T$ op\`ere sur $E_j(\tgamma)$ par le
caract\`ere $\twisted{w_j}{\mu}$.

Soit $i \in \{1,\dots, n\}$ et $j_1<j_2<\cdots<j_i$,
$j_{i+1}<\cdots<j_n$ la partition de $\{1,..,n\}$ telle que
$\mu(t_{j_s})=1$ pour $s \in \{1,\dots,i\}$ et
$\mu(t_{j_s})=-1$ pour $s \in \{i+1,\dots,n\}$.  Soit $\sigma
\in W_A$ la permutation telle que $\sigma*j_p=p$ pour $p \in
\{1,..,n\}$.  Observons que pour $p < q$,
$\sigma^{-1}(\alpha_{p,q}) = \alpha_{j_p,j_q} \in R_A^-$ si et
seulement si $j_p > j_q$, ce qui implique $1 \le p \le i < q
\le n$. Donc $R(\sigma^{-1})$ est contenu dans $\{\alpha_{p,q}
: 1 \le p \le i < q \le n\}$. Par cons\'equent
\begin{equation}
\label{eq2}
R(\sigma) = - \sigma^{-1}(R(\sigma^{-1})) \, \subset \,
\{\alpha_{j_q,j_p} : 1 \le p \le i < q \le n\}. 
\end{equation}
On pose
$$
\mu_i=\twisted{\sigma}{\mu}.
$$
Alors, $W_A(\mu_i)$ est le groupe de Weyl $W(R_i)$ associ\'e
au syst\`eme de racines $R_i$ de type $A_{i-1}\times
A_{n-i-1}$, de base
$S_i=\{\alpha_1,\dots,\alpha_{i-1},\alpha_{i+1},\dots,
\alpha_{n-1}\}$.  Donc $W_A(\mu)=\sigma W_A(\mu_i)
\sigma^{-1}=W_A(\sigma(R_i))$ admet pour base de racines
simples
$$
S(\mu) =\{ \alpha_{j_1,j_2}, \dots, \alpha_{j_{i-1},j_i} ,
\alpha_{j_{i+1},j_{i+2}}, \dots, \alpha_{j_{n-1},j_n}\}.
$$
On notera $\ell_\mu$ la longueur sur $W_A(\mu)$ associ\'ee
\`a $S(\mu)$.

\subsection{Restriction \`a $E_1(\tgamma)$}
\label{sec3.2}

Rappelons (voir~\ref{sec1.6}) que si $M \in \Rmod{\hb}$ et $w
\in W_A$, $\twisted{w}{M}$ est le $\H_B$-module tordu par
l'action de $\Int(w)$.
 
\begin{Prop}  
\label{prop4}  
Soit $w \in W_A$. On suppose que $\mu(\tk_\alpha)=0$ pour tout
$\alpha \in R(w^{-1})$.  Alors l'application $\phi :
M(\twisted{w^{-1}}{\tgamma}) \rightarrow
\twisted{w}{M(\tgamma)}$, donn\'ee par $\phi(t_g\otimes
v_{\twisted{w^{-1}}{\tgamma}})=t_{wgw^{-1}}\otimes v_\tgamma$
pour tout $g \in W_A$, est un isomorphisme de $\H_B$-modules.
\end{Prop}  
  
\begin{dem} 
  Comme l'espace vectoriel $M(\tgamma)=\twisted{w} M(\tgamma)$
  admet pour base $\{t_g \otimes v_{\tgamma}=t_{w^{-1}gw}
  \centerdot (1 \otimes v_{\tgamma}, g \in W\}$, le vecteur $1
  \otimes v_\tgamma$ engendre $\twisted{w}{M(\tgamma)}$.  Dans
  $\twisted{w}M(\tgamma)$ le groupe $T$ agit sur $1 \otimes
  v_\tgamma$ par le caractère $\twisted{w^{-1}}{\mu}$; de plus,
  pour tout $\zeta \in V$, il vient en utilisant le
  Lemme~\ref{lem1}:
  $$
  \zeta \centerdot 1 \otimes v_{\tgamma} = t_w \zeta
  t_{w^{-1}}. 1 \otimes v_{\tgamma}=\Bigl(w(\zeta) +
  \sum_{\alpha \in R(w^{-1})}<\alpha\spcheck,
  w(\zeta)>t_{s_{\alpha}} \tk_{\alpha}\Bigr).1 \otimes
  v_{\tgamma}= \twisted{w^{-1}}{\gamma}(\zeta).1 \otimes
  v_{\tgamma}
  $$
  puisque $\tk_{\alpha}.1 \otimes v_{\tgamma}=
  \mu(\tk_{\alpha}).1 \otimes v_{\tgamma}=0$ pour tout $\alpha
  \in R(w^{-1})$.  Donc $1 \otimes v_\tgamma$ est de poids
  $\twisted{w^{-1}}{\tgamma}$ dans $\twisted{w}{M(\tgamma)}$.
  Par la propri\'et\'e universelle de l'induction, il existe un
  unique $\H_B$-morphisme surjectif $\phi$ de
  $M(\twisted{w^{-1}}{\tgamma})$ dans
  $\twisted{w}{M(\tgamma)}$, qui envoie $1\otimes
  v_{\twisted{w^{-1}}{\tgamma}}$ sur $1 \otimes v_\tgamma$.
  Ces deux modules \'etant de dimension $|W_A|$, $\phi$ est
  bijectif.  Enfin, pour tout $g \in W_A$ on a
  $$
  \phi(t_g\otimes v_{\twisted{w^{-1}}{\tgamma}})=\phi(t_g.(1
  \otimes v_{\twisted{w^{-1}}{\tgamma}})) =t_g\centerdot \phi(1
  \otimes v_{\twisted{w^{-1}}{\tgamma}})=t_{wgw^{-1}}.(1
  \otimes v_{\tgamma})=t_{wgw^{-1}} \otimes v_{\tgamma},
  $$
  ce qui d\'emontre la proposition.
\end{dem}

\begin{cor}
\label{cor5}
{\rm (1)} On a $\mu(\tk_\alpha) =0$ pour tout $\alpha \in
R(\sigma)$ et l'application $\phi :
M(\twisted{\sigma}{\tgamma}) \rightarrow
\twisted{\sigma^{-1}}{M(\tgamma)}$ donn\'ee par
$\phi(t_g\otimes
v_{\twisted{\sigma}{\tgamma}})=t_{\sigma^{-1}g\sigma}\otimes
v_\tgamma$, pour tout $g \in W_A$, est un isomorphisme de
$\H_B$-modules.

\noindent {\rm (2)} Il existe un isomorphisme
$$
\twisted{\sigma}{\phi} : M(\twisted{\sigma}{\tgamma})
\longisomto M(\tgamma)
$$
tel que $\twisted{\sigma}{\phi}(t_g \otimes
v_{\twisted{\sigma}{\tgamma}}) = t_{g \sigma} \otimes
v_\tgamma$.
\end{cor}

\begin{dem}
  (1) Par~\eqref{eq2} on sait que toute racine $\alpha \in
  R(\sigma)$ s'\'ecrit $\alpha_{j_q,j_p}$ avec $1 \le p \le i <
  q \le n$. Par cons\'equent $\mu(\tk_\alpha)= 1
  +\mu(t_{j_p})\mu(t_{j_q})=0$. Il suffit donc d'appliquer la
  proposition pr\'ec\'edente \`a $w=\sigma^{-1}$.

\noindent (2) La Proposition~\ref{prop61} fournit un
isomorphisme de $\twisted{\sigma^{-1}}{M(\tgamma)}$ sur
$M(\tgamma)$ qui envoie $t_u \otimes v_\tgamma$ sur $t_{\sigma
  u} \otimes v_\tgamma$; en le composant avec $\phi$ on obtient
l'isomorphisme $\twisted{\sigma}{\phi}$ voulu.
\end{dem}

Le~(2) du Corollaire~\ref{cor5} montre que l'\'etude du module
$M(\tgamma)$ est \'equivalente \`a celle de
$M(\twisted{\sigma}{\tgamma})$. Nous verrons que $E_1(\tgamma)$
h\'erite d'une structure de $\H_B(\mu)$-module d\'eterminant
celle de $M(\tgamma)$ (cf.~Proposition~\ref{prop12}); il est
donc naturel de comparer les modules $E_1(\tgamma)$ et
$E_1(\twisted{\sigma}{\tgamma})$, ce que nous ferons au
Th\'eor\`eme~\ref{thm7}.

Observons que, puisque $W_A(\twisted{\sigma}{\mu}) = \sigma
W_A(\mu)\sigma^{-1}$, l'application $\phi$ du
Corollaire~\ref{cor5} donne un isomorphisme d'espaces
vectoriels
$$
\phi : E_1(\twisted{\sigma}{\tgamma}) \to E_1(\tgamma).
$$

\begin{Lemme}
\label{lem6}
La restriction de $\psi$ \`a $\H_B(\mu)$ munit $E_1(\tgamma)$
d'une structure de $\H_B(\mu)$-module.
\end{Lemme}

\begin{dem}
  L'espace $E_1(\tgamma)$ est \'evidemment stable sous l'action
  de $T$ et de $W_A(\mu)$. Il reste donc \`a montrer qu'il est
  stable sous l'action de $V$; puisque les \'el\'ements de $V$
  et $T$ commutent dans $\H_B$, ceci résulte du fait que
  $E_1(\tgamma)$ est la composante isotypique de type $\mu$ du
  $T$-module $M(\tilde{\gamma})$. Afin d'expliquer le rôle que
  va jouer l'algèbre de Hecke graduée $\H_{\mathit{gr}}(\mu_i)$
  au \S~\ref{sec3.3} nous donnons ci-dessous une autre preuve
  de ce fait.
  
  Soient $\zeta \in V$, $w \in W_A(\mu)$, et montrons par
  r\'ecurrence sur $\ell_\mu(w)$ que $\zeta.t_w \otimes
  v_{\tgamma} \in E_1(\tgamma)$.  En choisissant une
  d\'ecomposition r\'eduite de $w$ dans $W_A(\mu)$ on peut
  \'ecrire $w=s_{j_l,j_{l+1}}u$ avec $l \in \{1,\dots,n\}
  \smallsetminus \{i\}$ et $u \in W_A(\mu)$ tel que
  $\ell_\mu(u)< \ell_\mu(w)$.  Par d\'efinition des entiers
  $j_r$ on a, pour tout $j \in \{j_l+1,\dots, j_{l+1}-1\}$,
  $\mu(t_j)=-\mu(t_{j_l})=-\mu(t_{j_{l+1}})$.  Donc
  $\mu(\tk_{{j_l,j}}) = \mu(\tk_{{j,j_{l+1}}})=0$ et
  $\mu(\tk_{j_l,j_{l+1}}) = 2k$.  Comme
  $\mu(\tk_{u^{-1}(i),u^{-1}(j)})=\twisted{u}{\mu}(\tk_{i,j})
  =\mu(\tk_{i,j})$ et que $v_\tgamma$ est de poids $\mu$, le
  Lemme~\ref{lem1} appliqu\'e \`a $s_{p,q}=s_{j_l,j_{l+1}}$
  fournit:
  $$
  \zeta.t_w \otimes v_{\tgamma} =
  t_{s_{p,q}}.(s_{p,q}(\zeta).t_{u} \otimes
  v_{\tgamma})-2k<\alpha_{p,q}\spcheck,\zeta>t_{u} \otimes
  v_{\tgamma}.
  $$
  Par r\'ecurrence on a $s_{p,q}(\zeta).t_{u} \otimes
  v_{\tgamma} \in E_1(\tgamma)$, d'o\`u $\zeta.t_w \otimes
  v_{\tgamma} \in E_1(\tgamma)$.
\end{dem}

En appliquant le Lemme~\ref{lem6} avec
$\twisted{\sigma}{\tgamma}$ \`a la place de $\tgamma$, on
obtient que $E_1(\twisted{\sigma}{\tgamma})$ h\'erite
naturellement (via $\psi$) d'une structure de
$\H_B(\twisted{\sigma}{\mu})$-module.  Nous allons maintenant
transporter cette structure \`a $E_1(\tgamma)$ gr\^ace \`a
l'isomorphisme (d'espaces vectoriels) $\phi$ obtenu ci-dessus.

\begin{Theoreme}
\label{thm7}
{\rm (1)} Il existe un morphisme d'alg\`ebres
$$
\rho : \H_B(\twisted{\sigma}{\mu}) \rightarrow \End
E_1(\tgamma)
$$
d\'efini par $\rho(\xi)=\psi(\sigma^{-1}(\xi))$,
$\rho(t_w)=\psi(t_{\sigma^{-1}w\sigma})$, pour $\xi \in V$, $w
\in T \rtimes W_A(\twisted{\sigma}{\mu})$, qui donne une
structure de $\H_B(\twisted{\sigma}{\mu})$-module sur l'espace
vectoriel $E_1(\tgamma)$.  Pour cette structure, not\'ee
$E_1^{\sigma}(\tgamma)$, on a un isomorphisme de
$\H_B(\twisted{\sigma}{\mu})$-modules
$$
\vphi : E_1(\twisted{\sigma}{\tgamma}) \isomto
E_1^{\sigma}(\tgamma)
$$
donn\'e par $\vphi(t_g \otimes
v_{\twisted{\sigma}{\tgamma}}) =t_{\sigma^{-1}g\sigma} \otimes
v_\tgamma$, pour tout $g \in W_A(\twisted{\sigma}{\mu})$.
  
\noindent {\rm (2)} Un sous-espace $N \subset E_1(\tgamma)$ est un
sous-$\H_B(\mu)$-module si et seulement si c'est un
sous-$\H_B(\twisted{\sigma}{\mu})$-module de
$E_1^{\sigma}(\tgamma)$.
\end{Theoreme}

\begin{dem}
  (1) Nous allons montrer que le morphisme $\rho$ cherch\'e
  n'est autre que $\psi\circ\Int(\sigma^{-1})$ restreint \`a
  $\hb(\twisted{\sigma}{\mu})$. Observons tout d'abord que, par
  d\'efinition, $\psi(\Int(\sigma^{-1})(t_w)) =
  \psi(t_{\sigma^{-1}w\sigma}) \in \End E_1(\tgamma)$ pour tout
  $w \in T \rtimes W_A(\twisted{\sigma}{\mu})$. Il suffit donc
  de v\'erifier que $\psi(\Int(\sigma^{-1})(\xi)) =
  \psi(\sigma^{-1}(\xi))$ sur $E_1(\tgamma)$ pour tout $\xi \in
  V$. Pour ce faire, calculons $\psi(\Int(\sigma^{-1})(\xi)) =
  \psi(t_\sigma^{-1} \xi t_\sigma)$ sur un \'element $v \in
  E_1(\tgamma)$.  Comme dans la preuve de la
  Proposition~\ref{prop4} on obtient:
  $$
  t_\sigma^{-1} \zeta t_{\sigma}. v
  =\Bigl(\sigma^{-1}(\zeta) + \sum_{\alpha \in
    R(\sigma)}<\alpha\spcheck,
  \sigma^{-1}(\zeta)>t_{s_{\alpha}} \tk_{\alpha}\Bigr).v=
  \sigma^{-1}(\zeta).v
  $$
  puisque $\tk_{\alpha}.v= \mu(\tk_{\alpha}).v=0$ pour tout
  $\alpha \in R(\sigma)$ par le Corollaire~\ref{cor5} (on
  rappelle que $\CT$ op\`ere par $\mu$ sur $E_1(\tgamma)$).
  
  Il est clair que $\vphi$ est la bijection lin\'eaire $\phi :
  E_1(\twisted{\sigma}{\tgamma}) \to E_1(\tgamma)$
  pr\'ec\'edente. Il reste alors \`a v\'erifier que $\phi(a.y)=
  \rho(a)(\phi(y))$ pour tous $a \in
  \H_B(\twisted{\sigma}{\mu}), y \in
  E_1(\twisted{\sigma}{\tgamma})$. Mais, par le
  Corollaire~\ref{cor5}, on a $\phi(a.y)= a\centerdot \phi(y)=
  \psi(\Int(\sigma^{-1})(a))(y)= \rho(a)(y)$.

\noindent (2) Observons que $N$ est un
sous-$\H_B(\mu)$-module de $E_1(\tgamma)$ si, et seulement si,
$\{\xi.N \subset N, \, t_w.N \subset N\}$ pour tous $\xi \in
V$, $w \in T \rtimes W_A(\mu)$.  Ceci \'equivaut \`a
$\{\sigma^{-1}(\xi).N \subset N, \, t_{\sigma^{-1}(\sigma w
  \sigma^{-1})\sigma}.N \subset N\}$.  Compte tenu de (1),
cette condition s'\'ecrit encore $\{\rho(\xi)(N) \subset N, \,
\rho(t_{w'})(N) \subset N\}$ pour tous $\xi \in V$, $w' \in T
\rtimes W_A(\twisted{\sigma}{\mu})$ (puisque $T \rtimes
W_A(\twisted{\sigma}{\mu}) = \sigma(T \rtimes
W_A({\mu}))\sigma^{-1}$), ce qui signifie que $N$ est un
sous-$\H_B(\twisted{\sigma}{\mu})$-module de
$E_1^{\sigma}(\tgamma)$.
\end{dem}

\begin{cor}
\label{cor8}
On adopte les notations de~\ref{sec1.6} et du
Th\'eor\`eme~\ref{thm7}. Alors:
$$
\calL_{\H_B(\mu)}(E_1(\tgamma)) \cong
\calL_{\H_B(\twisted{\sigma}{\mu})}
(E_1(\twisted{\sigma}{\tgamma})), \qquad \lg_{\H_B(\mu)}
E_1(\tgamma) = \lg_{\H_B(\twisted{\sigma}{\mu})}
E_{1}(\twisted{\sigma}{\tgamma}).
$$
\end{cor}

\begin{dem}
  Par le (2) du Th\'eor\`eme~\ref{thm7} on sait que
  $\calL_{\H_B(\mu)}(E_1(\tgamma))=
  \calL_{\H_B(\twisted{\sigma}{\mu})}(E_1^\sigma({\tgamma}))$,
  et par le (1) du m\^eme th\'eor\`eme on obtient
  $\calL_{\H_B(\twisted{\sigma}{\mu})}
  (E_1(\twisted{\sigma}{\tgamma})) \cong
  \calL_{\H_B(\twisted{\sigma}{\mu})}(E_1^\sigma({\tgamma}))$.
  D'o\`u le corollaire.
\end{dem}

\subsection{Cas $\mu=\mu_i$}
\label{sec3.3}
On suppose dans cette section que $\tnu \in \calA$ est de la
forme $\nu\otimes\mu_i$.  On note $\H_{\mathit{gr}}(\mu_i)$
l'alg\`ebre de Hecke gradu\'ee associ\'ee au syst\`eme de
racines $R_i$ de $V$ et \`a la multiplicit\'e $c=2k$.  Pour
$\nu \in V_{\C}^*$, $M_i(\nu)$ d\'esignera le module de la
s\'erie principale sur $\H_{\mathit{gr}}(\mu_i)$ associ\'e \`a
$\nu$ (not\'e $M(\nu)$ en~\ref{sec1.1}).

Soit $I=\Ker(\mu_i)$ l'id\'eal de $\C T$ engendr\'e par les
$x-\mu_i(x)$, $x\in T$.  Comme les \'el\'ements de $V$
commutent \`a ceux de $T$, l'id\'eal \`a gauche $\H_B(\mu_i) I$
est bilat\`ere; on peut donc d\'efinir l'alg\`ebre quotient
$$
\hbi=\H_B(\mu_i) \otimes_{\C T} (\CT/I) =
\H_B(\mu_i)/\H_B(\mu_i)I.
$$
Le groupe $W_A(\mu_i)$ étant engendré par les réflexions
simples $s_{\alpha_p}$, $\alpha_p \in S_i$, il résulte
de~\eqref{eq1} et de (PBW3) que $\H_B(\mu_i)=S(V_\C)\otimes \C
W_A(\mu_i) \otimes \C T$.  Puisque $\C T/I$ est un $T$-module
de dimension~$1$ (de caract\`ere $\mu_i$), on a un isomorphisme
de $S(V_\C)$-modules $\hbi \simeq S(V_\C) \otimes \C
W_A(\mu_i)$.

Remarquons que pour tout $j \ne i$ on a
$\mu_i(t_j)=\mu_i(t_{j+1})$, donc $\mu_i(\tk_j)=2k$.  Les
relations de d\'efinition de $\H_B$ appliqu\'ees aux $\alpha_j
\in S_i$ fournissent alors dans $\hbi$:
\begin{equation*}
t_{s_{\alpha_j}}\zeta=s_{\alpha_j}(\zeta)t_{s_{\alpha_j}
}-2k<\alpha_j\spcheck,\zeta>, 
\end{equation*}
pour tout $\zeta \in V$.  Ces relations co\"{\i}ncident avec
celles d\'efinissant l'alg\`ebre de Hecke gradu\'ee
$\H_{\mathit{gr}}(\mu_i)$.  On en d\'eduit un morphisme
surjectif d'alg\`ebres
$$
F: \H_{\mathit{gr}}(\mu_i)\rightarrow \hbi
$$
d\'efini par $F(\zeta)=\zeta \pmod{\H_B(\mu_i)I}$ et
$F(t_w)=t_w \pmod{\H_B(\mu_i)I}$ pour $\zeta \in V$, $w \in
W_A(\mu_i)$. Puisque $\H_{\mathit{gr}}(\mu_i)$ et $\hbi$ sont
isomorphes \`a $S(V_\C) \otimes \C W_A(\mu_i)$ comme
$S(V_\C)$-modules, $F$ est un isomorphisme.

Observons que pour tout $\tnu=\nu \otimes \mu_i \in \calA$, $T$
op\`ere sur le $\H_B(\mu_i)$-module $E_1(\tnu)$ par le
caract\`ere $\mu_i$, donc $E_1(\tnu)$ peut \^etre consid\'er\'e
comme un $\hbi$-module.

\begin{Lemme}
\label{lem9}
L'application $f:M_i(\nu)\rightarrow E_1(\tnu)$ d\'efinie par
$f(t_w \otimes v_\nu)=t_w \otimes v_\tnu$ est un isomorphisme
qui entrelace, via $F$, les actions de
$\H_{\mathit{gr}}(\mu_i)$ et $\hbi$.
\end{Lemme}

\begin{dem}
  Gr\^ace \`a l'isomorphisme $F$, $E_1(\tnu)$ h\'erite d'une
  structure de $\H_{\mathit{gr}}(\mu_i)$-module pour laquelle
  il est engendr\'e par le vecteur $1 \otimes v_\tnu$, qui est
  de poids $\nu$ sous l'action de $S(V_\C)$.  La propri\'et\'e
  universelle de $M_i(\nu)$ assure donc l'existence et la
  surjectivit\'e de l'op\'erateur d'entrelacement $f$.  Comme
  $E_1(\tnu)$ et $M_i(\nu)$ sont de dimension $|W_A(\mu_i)|$,
  $f$ est bijective.
\end{dem}

On peut ainsi identifier le $\H_{\mathit{gr}}(\mu_i)$-module
$M_i(\nu)$ et le $\H_B(\mu_i)$-module $E_1(\tnu)$, il en
r\'esulte en particulier que:
\begin{equation}
\label{eq4}
\calL_{\H_{\mathit{gr}}(\mu_i)}(M_i(\nu)) \cong
\calL_{\H_B(\mu_i)}(E_1(\tnu)), \quad
\lg_{\H_{\mathit{gr}}(\mu_i)}M_i(\nu) =
\lg_{\H_B(\mu_i)} E_1(\tnu).
\end{equation}

\subsection{Crit\`ere d'irr\'eductibilit\'e}
\label{sec3.4}

Rappelons que $\sigma$ est tel que $\twisted{\sigma}{\mu}=
\mu_i$. Donc en appliquant~\eqref{eq4} \`a $\tnu =
\twisted{\sigma}{\tgamma} = \twisted{\sigma}{\gamma}\otimes
\mu_i$ et en utilisant le Corollaire~\ref{cor8} on obtient:

\begin{Prop}
 \label{prop11}
 M\^emes notations. On a
 $$
 \calL_{\H_B(\mu)}(E_1(\tgamma)) \cong
 \calL_{\H_{\mathit{gr}}(\mu_i)}(M_i(\twisted{\sigma}{\gamma})),
 \quad \lg_{\H_B(\mu)}E_1(\tgamma)=
 \lg_{\H_{\mathit{gr}}(\mu_i)}M_i(\twisted{\sigma}{\gamma}).
 $$
 En particulier, $E_1(\tgamma) $ est un $\hb(\mu)$-module
 simple si et seulement si $M_i(\twisted{\sigma}{\gamma})$ est
 un $\H_{\mathit{gr}}(\mu_i)$-module simple.
 \end{Prop}
 
 Rappelons que $\calS_R(M)$ désigne la suite des sous-quotients
 simples d'un $M \in \Rmod{R}$, cf.~\ref{sec1.6}, et que par
 \cite[Proposition2.8(c)]{KR}:
 $$
 \calS_{\H_{\mathit{gr}}(\mu_i)}(M_i(\twisted{\sigma}{\gamma}))
 =\calS_{\H_{\mathit{gr}}(\mu_i)}(M_i(\twisted{x\sigma}{\gamma}))
 $$
 pour tout $x \in W_A(\mu_i)$.

\begin{cor}
\label{cor11b}
Soit $w \in W_A$. Alors:
\begin{itemize}
\item il existe une bijection $f :
  \calS_{\H_B(\mu)}(E_1(\tgamma)) \to
  \calS_{\H_B(\twisted{w}{\mu})}(E_1(\twisted{w}{\tgamma}))$
  telle que $\dim f(Y) = \dim Y$ pour tout $Y \in
  \calS_{\H_B(\mu)}(E_1(\tgamma))$;
\item $E_1(\tgamma) $ est un $\hb(\mu)$-module simple si et
  seulement si $E_1(\twisted{w}{\tgamma})$ est un
  $\H_B(\twisted{w}{\mu})$-module simple.
\end{itemize}
\end{cor}

\begin{dem}
  Il existe $\varsigma \in W_A$ tel que $\twisted{\varsigma
    w}{\mu}=\mu_i$ et l'on peut écrire $\varsigma w = x \sigma$
  pour un $x \in W_A(\mu_i)$. Par la Proposition~\ref{prop11}
  on sait que $\calL_{\H_B(\mu)}(E_1(\tgamma)) \cong
  \calL_{\H_{\mathit{gr}}(\mu_i)}(M_i(\twisted{\sigma}{\gamma}))$
  et $\calL_{\H_B(\twisted{w}{\mu})}(E_1(\twisted{w}{\tgamma}))
  \cong \calL_{\H_{\mathit{gr}}(\mu_i)}(M_i(\twisted{\varsigma
    w}{\gamma}))$, l'existence de $f$ découle donc du rappel
  précédent. La deuxième assertion est alors évidente.
\end{dem}

Nous allons maintenant montrer que les treillis ordonnés
$\calL_{\hb}(M(\tgamma))$ et $\calL_{\hb(\mu)}(E_1(\tgamma))$
sont isomorphes.  Remarquons tout d'abord que (PBW3) et $\C
W_A(\mu) \otimes \CT \otimes S(V_\C) \subset \H_B(\mu)$
impliquent
\begin{equation*} 
  \H_B = \C W_A \otimes \CT \otimes S(V_\C) =
  \Bigl(\bigoplus_{j=1}^s t_{w_j} \C W_A(\mu)\Bigr) \otimes \CT
  \otimes
  S(V_\C)  \subset \sum_{j=1}^s t_{w_j}\H_B(\mu).
\end{equation*}
D'où:
\begin{equation} \label{eq5}
\H_B = \sum_{j=1}^s t_{w_j}\H_B(\mu).
\end{equation}

\begin{Prop}
\label{prop12}
Le module $M(\tgamma)$ s'identifie à $\H_B \otimes_{\H_B(\mu)}
E_1(\tgamma)$. Les applications $Y \to \hb \otimes_{\hb(\mu)}
Y$ et $X \to X \cap E_1(\tgamma)$ sont des bijections
réciproques de $\calL_{\hb(\mu)}(E_1(\tgamma))$ sur
$\calL_{\hb}(M(\tgamma))$.  On a en particulier:
$$
\lg_{\H_B} M({\tgamma}) =\lg_{\H_B(\mu)} E_1(\tgamma).
$$
\end{Prop}

\begin{dem}
  Soit $Y$ un sous-$\hb(\mu)$-module de $E_1(\tgamma)$.
  Puisque $M(\tgamma) = \bigoplus_{j=1}^s
  t_{w_j}.E_1(\tgamma)$, l'équation~\eqref{eq5} entraîne que le
  sous-module de $M(\tgamma)$ engendré par $Y$ est (comme
  sous-espace vectoriel) égal à
  $$
  \hb.Y = \sum_{j=1}^s t_{w_j}\H_B(\mu).Y =\sum_{j=1}^s
  t_{w_j}.Y= \bigoplus_{j=1}^st_{w_j}.Y.
  $$
  Il en résulte que l'application surjective canonique
  $$
  \pi : \hb \otimes_{\hb(\mu)} Y \sto \hb.Y, \quad \pi(a
  \otimes y) = a.y,
  $$
  est un isomorphisme. En effet, supposons $u =\sum_{l=1}^p
  a_l \otimes u_l \in \Ker \pi$. \'Ecrivons $a_l = \sum_{j=1}^s
  t_{w_j}a_{jl}$ avec $a_{jl} \in \hb(\mu)$.  Il vient $\pi(u)
  = \sum_{j,l} t_{w_j}a_{jl}.u_l = \sum_{j=1}^s
  t_{w_j}.\bigl(\sum_{l=1}^pa_{jl}.u_l \bigr)= 0$ avec
  $\sum_{l=1}^pa_{jl}.u_l \in Y$. Donc $\sum_{l=1}^p a_{jl}.u_l
  = 0$ pour tout $j$ et
  $$
  \textstyle{u=\sum_{j,l} t_{w_j} \otimes a_{jl}.u_l =
    \sum_jt_{w_j} \otimes \bigl(\sum_la_{jl}.u_l\bigr) = 0.}
  $$
  Nous pouvons ainsi identifier $\hb \otimes_{\hb(\mu)} Y$
  et $\hb.Y$.  Le foncteur $Y \to \hb \otimes_{\hb(\mu)} Y$
  fournit donc une application injective (qui préserve les
  inclusions) de $\calL_{\hb(\mu)}(E_1(\tgamma))$ dans
  $\calL_{\hb}(M(\tgamma))$.  Pour d\'emontrer la proposition
  il nous reste \`a voir que tout sous-module $X$ de
  $M(\tgamma)$ est de la forme $\hb \otimes_{\hb(\mu)} (X \cap
  E_1(\tgamma))$.  Soit $x \in X$, que l'on écrit $x
  =\sum_{j=1}^s t_{w_j}.x_j$ avec $x_j \in E_1(\tgamma)$ pour
  tout $j$. Par la Remarque~\ref{rem7} il existe des $y_i \in
  T$ tels que la matrice $[\twisted{w_j}{\mu}(t_{y_i})]_{i,j}$
  est inversible.  Puisque $T$ opère sur $t_{w_j}.x_j$ par le
  caractère $\twisted{w_j}{\mu}$, en appliquant les $t_{y_i}$ à
  $x$ on obtient que $t_{w_j}.x_j \in \CT.x \subset X$ pour
  tout $j$.  Il en résulte que $x_j \in X \cap E_1(\tgamma)$ et
  donc $X= \bigoplus_{j=1}^s t_{w_j}.(X \cap E_1(\tgamma))$,
  comme voulu.
\end{dem}

Observons que la proposition précédente assure que $M(\tgamma)$
est simple si et seulement si $E_1(\tgamma)$ l'est.  Nous
pouvons maintenant en déduire un critère d'irréductibilité pour
les modules $M(\tgamma)$. Rappelons que $\tgamma = \gamma
\otimes \mu \in \calA$ et que $\sigma \in W_A$ est a été
introduit en~\ref{sec3.1} de sorte que $\mu_i =
\twisted{\sigma}{\mu}$.  Comme dans~\cite[Theorem~2.10]{KR} on
pose pour tout $\nu \in V_\C^*$:
$$
P_i(\nu) =\{\alpha \in R_i^+ : \nu(\alpha)= \pm 2k\}.
$$

\begin{Theoreme}
\label{thm13}
Les assertions suivantes sont équivalentes:
\begin{enumerate}[{\rm (i)}]
\item $M(\tgamma)$ est un $\hb$-module simple;
\item $M(\twisted{w}{\tgamma})$ est un $\hb$-module simple pour
  tout $w \in W_A$;
\item $P_i(\twisted{\sigma}{\gamma}) = \emptyset$.
\end{enumerate}
\end{Theoreme}

\begin{dem}
  Par~\cite[Theorem~2.10]{KR} on sait que le
  $\H_{\mathit{gr}}(\mu_i)$-module
  $M_i(\twisted{\sigma}{\gamma})$ est simple si et seulement si
  $P_i(\twisted{\sigma}{\gamma}) = \emptyset$. Le th\'eor\`eme
  d\'ecoule donc de la combinaison de la
  Proposition~\ref{prop11}, du Corollaire~\ref{cor11b} et de la
  Proposition~\ref{prop12} (appliqu\'ee \`a $\tgamma$ et
  $\twisted{w}{\tgamma}$).
 \end{dem}

\section{Application à un système de racines de type $D_n$}
\label{sec4}

On d\'efinit dans cette section une sous-alg\`ebre $\hd$ de
$\H_B$ dont l'\'etude des repr\'esentations se d\'eduit du
travail pr\'ec\'edent. Elle peut se r\'ealiser dans l'alg\`ebre
de Cherednik rationnelle de type $D_n$: cette sous-alg\`ebre
est isomorphe à l'alg\`ebre $\D$ de l'introduction pour un
syst\`eme de racines de type $D_n$.

\subsection{L'algèbre $\hd$}
\label{sec4.1}
Rappelons que le système de racines $R_B$ contient le système
de racines $R_D$ suivant, de type $D_n$ et de base
$S_D=\{\alpha_i = \epsilon_i - \epsilon_{i+1}, 1 \le i \le n-1;
\, \epsilon_{n-1} + \epsilon_n\}$:
$$
R_D = \{\pm \epsilon_i \pm \epsilon_j : 1 \le i,j \le n\}
\supset R_D^+ = \{\epsilon_i \pm \epsilon_j : 1 \le i <j \le
n\}.
$$
Le groupe de Weyl $W_D$ associé est égal à $U \rtimes W_A$,
où $U$ est le sous-groupe de $T$ engendré par les
$s_{\epsilon_i}s_{\epsilon_{i+1}}$ pour $1 \le i \le n-1$.
Observons que, puisque $k \ne 0$, on a $\CU = \C[\tk_{\alpha_i}
: 1 \le i \le n-1]$.  Rappelons que l'automorphisme $\delta$
envoie $t_i$ sur $-t_i$; on a donc aussi $\CU = (\CT)^\delta =
\{a \in \CT : \delta(a) = a\}$. (Dans toute la suite on notera
$(\phantom{R})^\delta$ l'espace des $\delta$-invariants.)

On pose
\begin{equation*}
\hd = \hb^\delta = \{ h \in \hb : \delta(h) = h \}.
\end{equation*}
Comme $\hb = S(\VC) \otimes \C W_A \otimes \CT$ (par (PBW3)) et
que $\delta$ est l'identité sur $S(\VC)$ et $\C W_A$ on a
\begin{equation*}
\tag{PBW4}
\hd = S(\VC) \otimes \C W_A \otimes \CU.
\end{equation*}
D'autre part les relations (a), (b), (c) et
l'équation~(\ref{eq1}) du \S~\ref{sec1.2} montrent que
$$
\hd = \C \langle \xi \in V; \, w \in W_D \rangle = \C
\langle \xi \in V; \, t_{s_{\alpha_i}}, \, i=1,\dots,n-1
\rangle.
$$

Le calcul du centre de $\hd$ résulte de celui de $Z(\hb)$ fait
au \S~\ref{sec1.5}, dont on adopte les notations.

\begin{Theoreme}
\label{thm41}
On a: $Z(\hd) = S(\VC)^{W_A} \otimes (\CU)^{W_A} =
\bigoplus_{i=0}^{[\frac{n}{2}]}S(\VC)^{W_A} \vartheta_{2i}$.
\end{Theoreme}

\begin{dem}
  Observons tout d'abord que $\delta(\vartheta_j)=(-1)^j
  \vartheta_j$, donc
  $$
  \textstyle{ (\CU)^{W_A} = \CU \cap (\CT)^{W_A}
    =(\CT)^\delta \cap (\CT)^{W_A} =
    \bigoplus_{i=0}^{[\frac{n}{2}]} \C \vartheta_{2i}.  }
  $$
  Soit $c \in Z(\hd)$, c'est à dire $c \in \hd$ et $[c,V] =
  [t_{s_{\alpha_i}},c] = 0$ pour $i=1,\dots,n-1$. Par le
  Lemme~\ref{lem3} la condition $[c,V] =0$ équivaut à $c \in
  (S(\VC) \otimes \CT)^\delta = S(\VC) \otimes \CU$. La preuve
  du Théorème~\ref{thm4} montre alors que $c \in Z(\hd)$ si et
  seulement si $c \in Z(\hb)^\delta = \bigl( \bigoplus_{j=0}^n
  S(\VC)^{W_A} \vartheta_j\bigr)^\delta =
  \bigoplus_{i=0}^{[\frac{n}{2}]}S(\VC)^{W_A} \vartheta_{2i}$.
\end{dem}

\subsection{Réalisation dans l'algèbre de Cherednik}
\label{sec4.2}

Rappelons (Théorème~\ref{PBW3}) l'existence de l'isomorphisme
$\Phi : \hb \to \D_B$ tel que $\Phi(\epsilon_j)= D_j=D^B_j$ et
$\Phi(w) = w$ pour $w \in W_B$. Par restriction $\Phi$ induit
donc un isomorphisme de $\hd$ sur la sous-algèbre
$$
\D_D = \C\langle D_1^B,\dots, D_n^B; w \in W_D \rangle = \C
\langle D_1^B,\dots, D_n^B; w \in W_A \rangle.
$$
On peut ainsi réaliser $\D_D$ comme sous-algèbre de
l'algèbre de Cherednik $\calH_B(\bk)$ de type $B_n$, où $\bk=
(k,k_c)$ comme au \S~\ref{sec1.3}.  Notons que les relations de
d\'efinition de $\H_B$ ne font intervenir la multiplicit\'e
$\bk$ que par sa valeur $k$ sur les racines longues. Ceci
permet de réaliser $\hd$ comme sous-algèbre de l'algèbre de
Cherednik rationnelle $\calH_D(k)$ de type $D_n$. En effet,
choisissons pour multiplicité $\bk_0=(k,0)$ (i.e.~$\bk$ nulle
sur les racines courtes). L'opérateur de Dunkl $T_y^B(\bk_0)$,
pour $y \in \fa$, est alors
$$
\textstyle{ T_y^D(k) = \partial_y + \frac{k}{2} \sum_{\alpha
    \in R_D}\frac{<{\alpha},{y}>}{\alpha} (1 - t_{s_\alpha}), }
$$
qui est égal à l'opérateur de Dunkl défini par $y$ pour le
type $D_n$. L'algèbre de Cherednik $\calH_D(k)$ étant
engendrée, dans $\End_\C S(\fa^*)$, par les $T_y^D(k)$, $y \in
\fa$, $x \in \fa^*$ et $t_w$, $w \in W_D$, s'identifie donc à
une sous-algèbre de $\calH_B(\bk_0)$. Il est clair que l'image
$\Phi(\hd)$ de l'isomorphisme précédent est alors contenue dans
$\calH_D(k)$ et co\"{\i}ncide avec la sous-algèbre $\D$ définie
dans l'introduction.

\subsection{Modules de la s\'erie principale}
\label{sec4.3}

Notons $\mu \mapsto \bar{\mu}$ le morphisme (surjectif) de
restriction entre les groupes de caractères $T\spcheck$ et
$U\spcheck$. Les antécédents de $\bmu \in U\spcheck$ sont donc
$\mu$ et $-\mu = \sgn \otimes \mu$. On désigne par $W_A(\bmu) =
\{ w \in W_A : \twisted{w}{\bmu} = \bmu \}$ le stabilisateur de
$\bmu$ dans $W_A$. Remarquons que $\twisted{w}{\bmu} =
\overline{{\twisted{w}{\mu}}}$ pour tout $w \in W$ et que
$W_A(\mu) \subset W_A(\bmu)$. On définit deux sous-algèbres de
$\hd$ par
$$
\hd(\bmu) = \C \langle V, \, W_A(\bmu), \, U \rangle \supset
\hd(\mu) = \C \langle V, \, W_A(\mu), \, U \rangle.
$$
Observons que le noyau $\bar{I}$ du caractère $\bmu$
engendre un idéal bilatère de $\hd(\bmu)$ et que l'on peut
ainsi définir les algèbres quotients
$$
\overline{\hd(\bmu)} = \hd(\bmu) /\hd(\bmu) \bar{I}, \quad
\overline{\hd(\mu)} = \hd(\mu) /\hd(\mu)\bar{I}.
$$

Fixons $\mu \in T\spcheck$ et $\sigma \in W_A$ tel que
$\twisted{\sigma}{\mu} = \mu_i$ comme au \S~\ref{sec3.1}. On a
donc $W_A({\mu_i}) = \sigma W_A(\mu)\sigma^{-1}$ et
$W_A({\bmu_i}) = \sigma W_A(\bmu)\sigma^{-1}$. Rappelons que
$W_A(\mu_i)$ est le groupe de Weyl associ\'e au syst\`eme de
racines $R_i$ de type $A_{i-1} \times A_{n-i-1}$ de base $S_i
=\{\alpha_1,\dots,\alpha_{i-1},\alpha_{i+1},\dots,
\alpha_{n-1}\}$.

\begin{lem}
\label{lem42}
{\rm (1)} On a $W_A(\mu) \ne W_A(\bmu)$ si et seulement si
$n=2m$ et $\twisted{\sigma}{\mu} = \mu_m$. Dans ce cas
$W_A(\bmu_m) = W_A(\mu_m) \rtimes \langle w_0 \rangle$ où $w_0$
est l'élément le plus long de $W_A$.

\noindent {\rm (2)} L'application
$F: \H_{\mathit{gr}}(\mu_i)\rightarrow \overline{\hd(\mu_i)}$
d\'efinie par $F(\zeta)=\zeta \pmod{\H_D(\mu_i)\bar{I}}$ et
$F(t_w)=t_w \pmod{\H_D(\mu_i)\bar{I}}$ pour $\zeta \in V$, $w
\in W_A(\mu_i)$ est un isomorphisme d'algèbres. On a en
particulier:
$$
\overline{\hd(\mu_i)} \simeq \overline{\hb({\mu_i})} \simeq
\H_{\mathit{gr}}(\mu_i).
$$
\end{lem}

\begin{dem}
  (1) On a vu que $W_A(\mu) \ne W_A(\bmu)$ équivaut à
  $W_A(\mu_i) \ne W_A(\bmu_i)$.  Soit $w \in W_A(\bmu_i)
  \smallsetminus W_A(\mu_i)$, c'est à dire $\twisted{w}{\mu_i}=
  -\mu_i$. De $\twisted{w}{\mu_i}(t_j) =\mu_i(t_{w^{-1}*j})= -
  \mu(t_j)$ pour tout $j$ on tire que $\sum_j \mu_i(t_j)=
  \sum_j \mu_i(t_{w^{-1}*j}) = -\sum_j \mu_i(t_j)$ et donc
  $\sum_j \mu_i(t_j) = 0$. Cette condition signifie clairement
  que $n=2m$ et $i=m$. Inversement si cette condition est
  vérifiée l'élément $w_0$, qui envoie $t_j$ sur $t_{n+1 - j}$,
  vérifie $\twisted{w}{\mu_m}= -\mu_m$ et par conséquent $w_0
  \in W_A(\bmu_m) \smallsetminus W_A(\mu_m)$.  L'égalité
  $W_A(\bmu_m) = W_A(\mu_m) \rtimes \langle w_0 \rangle$ est
  alors triviale.
  
  \noindent (2) On procède comme au \S~\ref{sec3.3}. On observe
  que le groupe $W_A(\mu_i)$ est engendré par les réflexions
  simples $s_{\alpha_p}$, $\alpha_p \in S_i$, et que par (PBW4)
  il en résulte $\H_D(\mu_i)=S(V_\C)\otimes \C W_A(\mu_i)
  \otimes \C U$.  Puisque $\C U/\bar{I}$ est de dimension~$1$
  on a un isomorphisme de $S(V_\C)$-modules
  $\overline{\hd(\mu_i)} \simeq S(V_\C) \otimes \C W_A(\mu_i)$.
  D'autre part, les relations de d\'efinition de $\H_D$
  appliqu\'ees aux $\alpha_j \in S_i$ étant les mêmes que
  celles définissant $\H_{\mathit{gr}}(\mu_i)$, l'application
  $F$ est bien définie et est un isomorphisme. Le dernier
  isomorphisme a \'et\'e vu au \S~\ref{sec3.3}.
\end{dem}

Notons $\bar{\calA}=\{ \gamma\otimes\bnu : \gamma \in V_\C^*,
\, \bnu \in U\spcheck\}$ l'ensemble des morphismes d'alg\`ebres
de $S(V_\C)\otimes\CU$ vers $\C$.  Un élément $w \in W_A$ opère
sur $\bar{\calA}$ par $\twisted{w}(\gamma\otimes\bnu) =
\twisted{w}{\gamma} \otimes \twisted{w}{\bnu}$.  Pour un
$\hd$-module de dimension finie $M$ on d\'efinit de m\^eme
qu'au \S~\ref{sec2.1} les sous-espaces poids et sous-espaces
poids g\'en\'eralis\'es associ\'es \`a un \'el\'ement de
$\bar{\calA}$.  Comme $S(V_\C)\otimes\CU$ est une
sous-alg\`ebre commutative de $\hd$, $M$ est somme directe de
ses sous-espaces poids g\'en\'eralis\'es.
 
Fixons $\bgamma=\gamma\otimes \bmu \in \bar{\calA}$ et
d\'efinissons maintenant le $\hd$-module de la s\'erie
principale $N(\bgamma)$ associ\'e. Si $\C v_\bgamma$ est le
$S(V_\C)\otimes\CU$-module de dimension $1$ défini par
$\bgamma$, on induit ce module à $\hd$:
$$
N(\bgamma) =\hd \otimes_{S(V_\C)\otimes\C U} \C v_\bgamma.
$$
Ce module admet pour base $\{t_w \otimes v_{\bgamma} : w \in
W_{A}\}$, ordonn\'ee comme pr\'ec\'edemment. Dans cette base
les éléments de $V_\C \bigoplus \CU$ admettent une matrice
triangulaire sup\'erieure et les poids de $N(\bgamma)$ sont
donc les $\twisted{w}{\bgamma}$, $w \in W_A\}$.  On notera
encore $\psi : \hd \to \End_\C N(\bgamma)$ la repr\'esentation
d\'efinie par $N(\bgamma)$.

L'\'etude des $\hd$-modules simples se ram\`ene \`a l'\'etude
des $\hd$-modules de la s\'erie principale car on montre de
m\^eme que pour $\H_B$:

\begin{Prop}
\label{prop43}
Soit $M$ un $\hd$-module simple et $\bgamma \in \bar{\calA}$ un
poids de $M$.  Alors $M$ est quotient de $N(\bgamma)$. En
particulier tout $\hd$-module simple est de dimension au plus
$|W_A|$.
\end{Prop}

Posons $\tgamma = \gamma \otimes \mu \in \calA$. Similairement
à ce qui a été fait pour $M(\tgamma)$ on peut d\'ecomposer
$N(\bgamma)$ en composantes isotypiques sous l'action de $U$.
Soient $w_1=\id,\dots,w_r$ des repr\'esentants des classes \`a
gauche de $W_A$ modulo $W_A(\bmu)$; posons
$$
F_j(\bgamma) = \bigoplus_{w \in W_A(\bmu)} \C t_{w_jw}
\otimes v_\bgamma.
$$
On a alors
$$
N(\bgamma) = \bigoplus_{j=1}^r F_j(\bgamma), \quad
F_j(\bgamma) = t_{w_j}. F_1(\bgamma),
$$
et le groupe $U$ op\`ere sur $F_j(\bgamma)$ par le
caract\`ere $\twisted{w_j}{\bmu}$. Observons aussi que les
$\twisted{w_j}{\bmu}$, $j=1,\dots,r$, forment une famille libre
dans le dual $(\CU)^*$ (de m\^eme que les $\twisted{w_j}{\mu}
\in (\CT)^*$ de la Remarque~\ref{rem7}).  Puisque les
\'el\'ements de $S(\VC)$ et $U$ commutent on obtient comme pour
$E_1(\tgamma)$ (voir Lemme~\ref{lem6}):

\begin{lem}
\label{lem43}
La restriction de $\psi$ \`a $\hd(\bmu)$ munit $F_1(\bgamma)$
d'une structure de $\hd(\bmu)$-module.
\end{lem}

Comme $\hd \subset \hb$, le $\hb$-module $M(\tgamma)$ est muni
par restriction d'une structure de $\hd$-module.

\begin{prop}
\label{prop44}
{\rm (1)} Le $\hd$-module $M(\tgamma)$ s'identifie \`a
$N(\bgamma)$ par $t_w \otimes v_\tgamma \mapsto t_w \otimes
v_\bgamma$.

\noindent {\rm (2)} Soit $w \in W_A$; les $\hd$-modules
$N(\bgamma)$ et $\twisted{w}{N(\bgamma)}$ sont isomorphes.

\noindent {\rm (3)}  Il existe un isomorphisme
$$
\twisted{\sigma}{\phi} : N(\twisted{\sigma}{\bgamma})
\longisomto N(\bgamma)
$$
tel que $\twisted{\sigma}{\phi}(t_g \otimes
v_{\twisted{\sigma}{\bgamma}}) = t_{g \sigma} \otimes
v_\bgamma$.

\noindent {\rm (4)} Le module $N(\bgamma)$ s'identifie à $\H_D
\otimes_{\H_B(\bmu)} F_1(\bgamma)$. Les applications $Y \to \hd
\otimes_{\hd(\bmu)} Y$ et $X \to X \cap F_1(\bgamma)$ sont des
bijections réciproques de $\calL_{\hd(\bmu)}(F_1(\bgamma))$ sur
$\calL_{\hd}(N(\bgamma))$.

\noindent {\rm (5)} Le $\hd$-module $N(\bgamma)$ est simple si
et seulement si le $\H_D(\bmu_i)$-module
$F_1(\twisted{\sigma}{\bgamma})$ l'est.
\end{prop}

\begin{dem}
  (1) De $M(\tgamma) = \bigoplus_{w \in W_A} t_w \otimes
  v_\tgamma$ on tire que $1\otimes v_\tgamma$ engendre
  $M(\tgamma)$ comme $\hd$-module. D'autre part $T$ op\'erant
  sur $1\otimes v_\tgamma$ via le caract\`ere $\mu$, le groupe
  $U$ agit sur cet \'el\'ement via $\bmu$. Donc $1\otimes
  v_\tgamma$ est de poids $\bgamma$ dans le $\hd$-module
  $M(\tgamma)$. La propri\'et\'e universelle de l'induction et
  le fait que $\dim M(\tgamma) = \dim N(\bgamma)$ assurent
  alors de l'existence d'un isomorphisme de $N(\bgamma)$ sur
  $M(\tgamma)$ tel que $t_w \otimes v_\bgamma \mapsto t_w
  \otimes t_\tgamma$ pour tout $w \in W_A$.

\noindent (2) Il suffit d'appliquer la
Proposition~\ref{prop61} en utilisant l'isomorphisme d\'ecrit
en~(1).

\noindent (3) r\'esulte de (1) et du Corollaire~\ref{cor5}. 

\noindent (4) La preuve est la m\^eme que celle de la
Proposition~\ref{prop12} en rempla\c{c}ant $\hb$ par $\hd$, $T$
par $U$, $W_A(\mu)$ par $W_A(\bmu)$ et en utilisant
l'ind\'ependance lin\'eaire des $\twisted{w_j}{\bmu}$,
$j=1,\dots,r$.

\noindent (5) d\'ecoule de (3) et (4). 
\end{dem}

Remarquons que $F_1(\bgamma) = \bigoplus_{w \in W_A(\bmu)} t_w
\otimes v_\bgamma$ est un $\overline{\H_D(\bmu)}$-module
contenant le sous-espace $\bigoplus_{w \in W_A(\mu)} \C \,t_w
\otimes v_\bgamma$. Ce dernier est un
$\overline{\hd(\mu)}$-module qui s'identifie
(cf.~Proposition~\ref{prop44}(1)) au sous-espace
$E_1(\tgamma)$.  Comme $\psi(\hd(\mu)) = \psi(\hb(\mu)) \subset
\End_\C E_1(\tgamma)$ (car ces alg\`ebres sont toutes deux
engendr\'ees par $\psi(V)$ et $\psi(W_A(\mu))$) on peut
indiff\'eremment consid\'erer $E_1(\tgamma)$ comme un
$\overline{\hd(\mu)}$ ou $\overline{\hb(\mu)}$-module.

\subsection{Crit\`ere d'irr\'eductibilit\'e de $F_1(\bgamma)$}  
\label{sec4.4}
Le (5) de la Proposition~\ref{prop44} ram\`ene la question de
l'irr\'eductibilit\'e du module $N(\bgamma)$ \`a celle du
$\H_D(\bmu_i)$-module $F_1(\twisted{\sigma}{\bgamma})$ avec
$\twisted{\sigma}{\bgamma}= \twisted{\sigma}{\gamma} \otimes
\mu_i$.  On suppose donc dans ce paragraphe que $\bgamma =
\gamma \otimes \mu$ et $\tgamma = \gamma \otimes \bmu$ avec
$\mu = \mu_i$. Rappelons (Lemme~\ref{lem42}) qu'on a alors
$$
\overline{\hd(\mu)} \simeq \overline{\hb({\mu})} \simeq
\H_{\mathit{gr}}(\mu),
$$
où $\H_{\mathit{gr}}(\mu)$ est l'algèbre de Hecke graduée
$\H_{\mathit{gr}}(2k,S_i)$, et que le $\hb(\mu)$-module
$E_1(\tgamma)$ s'identifie au $\H_{\mathit{gr}}(\mu)$-module de
la série principale $M(\gamma)$ (cf.~Lemme~\ref{lem9}).
\subsubsection{Le cas $W_A(\mu)=W_A(\bmu)$} 
\label{sec4.41}
On a ici $\hd(\bmu)=\hd(\bmu)$ et donc $F_1(\bgamma) =
E_1(\tgamma)$ comme $\H_{\mathit{gr}}(\mu)$-module. Par
conséquent:
\begin{equation}
\label{eq6}
\text{$F_1(\bgamma)$ simple $\iff$ $E_1(\tgamma)$ simple $\iff$
$M(\gamma)$ simple $\iff$ $\calP_i(\gamma) = \emptyset$.}
\end{equation}

\subsubsection{Le cas $W_A(\mu)\ne W_A(\bmu)$} 
\label{sec4.42}
Par le Lemme~\ref{lem42} on sait que: $n=2m$, $\mu = \mu_m$,
$W_A(\bmu) = W_A(\mu) \rtimes \langle w_0 \rangle$ où $w_0$ est
le plus grand élément de $W_A$. Le système de racines $R_m$ est
de type $A_{m-1} \times A_{m-1}$ et a pour groupe de Weyl
$W_m=W(R_m) \simeq W_A(\mu)\simeq \mathfrak{S}_m \times
\mathfrak{S}_m$. On note $\varpi_0$ le plus grand élément de
$W_m$; on a donc
$$
\varpi_0(\epsilon_j) = \epsilon_{m+1 - j}, \quad
\varpi_0(\epsilon_{m+j}) = \epsilon_{2m+1 - j}
$$
pour $j=1,\dots,m$. Définissons une involution $\tau \in
\GL(V)$ par
$$
\tau(\epsilon_j) = \epsilon_{m+j}, \, \text{pour tout
  $j\in\{1,\dots, m\}$.}
$$
Comme $\tau(\alpha_{p,q}) = \alpha_{p+m,q+m}$ pour $1 \le p
< q \le m$, l'application $\tau$ permute les éléments de
$R_m^+$ et induit un automorphisme (extérieur), $w \mapsto
\tau(w)=\tau w \tau$, de $W_m$ qui vérifie
$$
< \tau(\alpha)\spcheck, \tau(\xi) > = < \alpha\spcheck, \xi
>, \quad \tau(s_\alpha) = \tau s_\alpha \tau =
s_{\tau(\alpha)},
$$
pour tous $\alpha \in R_m$, $\xi \in V$. Remarquons aussi
que l'automorphisme $\omega_0: w \mapsto w_0 w w_0$ de $W_m$
s'écrit
$$
\omega_0 = \Int(\varpi_0) \circ \tau = \tau \circ
\Int(\varpi_0).
$$

\begin{Assertion}
\label{ass45}
{\rm (1)} L'automorphisme $\omega_0$ induit un automorphisme
involutif de $\H_{\mathit{gr}}(2k,S_m)=\H_{\mathit{gr}}(\mu)$
encore noté $\omega_0$ et donné par
$$
\omega_0(\xi) = w_0(\xi) + 2k \sum_{\alpha \in R_m^+} <
\alpha\spcheck, w_0(\xi) > t_{s_\alpha}, \quad \omega_0(t_w)=
t_{\omega_0(w)},
$$
pour tous $\xi \in V$ et $w \in W_m$.

\noindent {\rm (2)} Le module $\twisted{\omega_0}{M(\gamma)}$
est isomorphe à $M(\twisted{\tau}{\gamma})$, où
$\twisted{\tau}{\gamma} \in \VC^*$ est défini par
$\twisted{\tau}{\gamma}(\xi) = \gamma(\tau(\xi))$.
\end{Assertion}

\begin{dem}
  (1) Commen\c{c}ons par étendre $\tau \in \GL(V)$ en un
  automorphisme involutif de $\H_{\mathit{gr}}(\mu)$ en posant
  $\tau(t_w)= t_{\tau(w)}$ pour $w \in W_m$. Il suffit pour
  cela de vérifier les relations de définition de
  $\H_{\mathit{gr}}(2k,S_m)$ sont respect\'ees, ce qui résulte
  d'un calcul facile.  Observons ensuite que $\Int(\varpi_0) :
  a \mapsto t_{\varpi_0}a t_{\varpi_0}$ est un automorphisme
  involutif de $\H_{\mathit{gr}}(\mu)$ qui est donné par
  $$
  \Int(\varpi_0)(\xi) = \varpi_0(\xi) + 2k \sum_{\alpha \in
    R_m^+} < \alpha\spcheck, \varpi_0(\xi) > t_{s_\alpha},
  \quad \Int(\varpi_0)(t_w)= t_{\varpi_0 w \varpi_0},
  $$
  pour tous $\xi \in V$ et $w \in W_m$ (voir le
  Lemme~\ref{lem1}). Il reste alors à remarquer que
  $$
  \Int(\varpi_0) \circ \tau = \tau \circ \Int(\varpi_0) \in
  \Aut(\H_{\mathit{gr}}(\mu))
  $$
  est l'automorphisme $\omega_0$ cherché.

\noindent (2) La définition de $\twisted{\omega_0}{M(\gamma)}$
et l'égalité $\omega_0 = \tau \circ \Int(\varpi_0)$ montrent
que ce module est isomorphe à
$\twisted{\varpi_0}{(\twisted{\tau}{M(\gamma)})}$. Si l'on
prouve que $\twisted{\tau}{M(\gamma)} \simeq
M(\twisted{\tau}{\gamma})$ la Proposition~\ref{prop62} donnera
alors l'isomorphisme voulu. Observons que $1 \otimes v_\gamma
\in \twisted{\tau}{M(\gamma)}$ vérifie $\xi\centerdot 1 \otimes
v_\gamma = \twisted{\tau}{\gamma}(\xi).1 \otimes v_\gamma$ pour
$\xi \in \VC$ et est donc de poids $\twisted{\tau}{\gamma}$. Il
en découle l'existence d'un morphisme surjectif
$M(\twisted{\tau}{\gamma}) \sto \twisted{\tau}{M(\gamma)}$, qui
est un isomorphisme puisque ces deux modules sont de même
dimension.
\end{dem}

Définissons deux sous-algèbres de $\End_\C F_1(\bgamma)$ par
$$
\mbH(\bmu) = \psi(\hd(\bmu)) \supset \mbH(\mu) =
\psi(\hd(\mu)).
$$
L'isomorphisme $\H_{\mathit{gr}}(\mu) \isomto
\overline{\hd(\mu)}$ induit un morphisme surjectif d'algèbres
$\psi : \H_{\mathit{gr}}(\mu)\sto \mbH(\mu)$ tel que $\xi
\mapsto \psi(\xi)$ et $t_w \mapsto \psi(t_w)$ pour $\xi \in V,
w \in W_A(\mu)$. Dans $\hd$ l'automorphisme $\Int(w_0)$ est
défini par $\Int(w_0)(t_w) = t_{w_0 w w_0}$ et, pour $\xi \in
V$,
$$
\Int(w_0)(\xi)=w_0(\xi) + \sum_{\alpha \in R_D^+} <
\alpha\spcheck, w_0(\xi) > t_{s_\alpha}\tk_{\alpha}.
$$
Remarquons que $\psi(\tk_{\alpha})_{\mid F_1(\bgamma)} = 0$
lorsque $\alpha \notin R_m$, car $\bmu(t_pt_{m+q}) = -1$ pour
$1 \le p,q \le m$, et $\psi(\tk_{\alpha})_{\mid F_1(\bgamma)} =
2k.\id_{F_1(\bgamma)}$ si $\alpha \in R_m$. Par conséquent
$\Int(w_0)$ donne un automorphisme de $\mbH(\bmu)$, $\Int(w_0)
: a \mapsto \psi(t_{w_0})a \psi(t_{w_0})$ vérifiant
$$
\Int(w_0)(\psi(\xi))=\psi(w_0(\xi)) + 2k \sum_{\alpha \in
  R_m^+} < \alpha\spcheck, w_0(\xi) > \psi(t_{s_\alpha}).
$$
Il découle alors de l'Assertion~\ref{ass45}(1) que $\psi :
\H_{\mathit{gr}}(\mu)\sto \mbH(\mu)$ est tel que
$$
\psi(\omega_0(a)) = \Int(w_0)(\psi(a)) \ \, \text{pour tout
  $a \in \H_{\mathit{gr}}(\mu)$}.
$$

Comme $F_1(\bgamma)$ est un $\mbH(\bmu)$-module, nous pouvons
le considérer comme un $\H_{\mathit{gr}}(\mu)$-module grâce à
$\psi :\H_{\mathit{gr}}(\mu)\sto \mbH(\mu) \subset \mbH(\bmu)$.
Nous allons maintenant étudier ce
$\H_{\mathit{gr}}(\mu)$-module.  Via l'identification de
$E_1(\tgamma)$ avec $\bigoplus_{w \in W_A(\mu)} \C t_w \otimes
v_\bgamma$ on a
$$
F_1(\bgamma) = E_1(\tgamma) \bigoplus t_{w_0}. E_1(\tgamma).
$$
Soit $a \in \H_{\mathit{gr}}(\mu)$ et $x \in E_1(\tgamma)$;
il vient $\psi(a).(t_{w_0}.x) = t_{w_0}.(\Int(w_0)(\psi(a)))$.
Ceci montre que $t_{w_0}. E_1(\tgamma)$ est un
$\mbH(\mu)$-module isomorphe au module tordu
$\twisted{\Int(w_0)}{E_1(\tgamma)}$. Compte tenu de ce qui
précède, ce dernier module, regardé comme
$\H_{\mathit{gr}}(\mu)$-module, est isomorphe à
$\twisted{\omega_0}{M(\gamma)}$. Par l'Assertion~\ref{ass45}(2)
on a ainsi prouvé:

\begin{Assertion}
\label{ass46}
Le $\H_{\mathit{gr}}(\mu)$-module $F_1(\tgamma)$ est isomorphe
à $M(\gamma) \bigoplus M(\twisted{\tau}{\gamma})$.
\end{Assertion}

Lorsque $M(\gamma)$ est simple la nature du $\hd(\bmu)$-module
$F_1(\bgamma)$ résulte de l'assertion qui suit.

\begin{Assertion}
\label{ass47}
Supposons le $\H_{\mathit{gr}}(2k,S_m)$-module $M(\gamma)$
irréductible. Alors, le $\hd(\bmu)$-module $F_1(\bgamma)$ n'est
pas simple si et seulement si $M(\gamma) \simeq
M(\twisted{\tau}{\gamma})$.
\end{Assertion}

\begin{dem}
  Observons tout d'abord que l'hypothèse et
  l'Assertion~\ref{ass46} assurent que $F_1(\bgamma)$ est un
  $\H_{\mathit{gr}}(\mu)$-module semi-simple de longueur $2$,
  isomorphe à $M(\gamma) \bigoplus M(\twisted{\tau}{\gamma})$.
  
  Soit $V$ un sous-$\hd(\bmu)$-module de $F_1(\bgamma)$ non nul
  et différent de $F_1(\bgamma)$. En tant que
  $\H_{\mathit{gr}}(\mu)$-module on a par conséquent $\lg V =1$
  et l'on est dans l'un des deux cas suivants:
\begin{itemize}
\item $F_1(\bgamma)= V \bigoplus t_{w_0}.E_1(\tgamma)$ et donc
  $V \simeq M(\gamma)$;
\item $F_1(\bgamma)= V \bigoplus E_1(\tgamma)$ et donc $V
  \simeq M(\twisted{\tau}{\gamma})$.
\end{itemize}
Pla\c{c}ons nous dans le premier cas (le deuxième est
similaire). Puisque $t_{w_0}.F_1(\bgamma) = F_1(\bgamma)$ et
$t_{w_0}.V = V$ on obtient $F_1(\bgamma) = V \bigoplus
t_{w_0}.E_1(\tgamma)$ et donc $V \simeq t_{w_0}.E_1(\tgamma)
\simeq M(\twisted{\tau}{\gamma})$ comme
$\H_{\mathit{gr}}(\mu)$-module. D'où $M(\gamma) \simeq
M(\twisted{\tau}{\gamma})$.

Réciproquement supposons $M(\gamma) \simeq
M(\twisted{\tau}{\gamma})$. On a donc $E_1(\tgamma) \simeq
\twisted{\Int(w_0)}{E_1(\tgamma)}$, simple comme
$\mbH(\mu)$-module. Fixons un isomorphisme $f : E_1(\tgamma)
\isomto \twisted{\Int(w_0)}{E_1(\tgamma)}$. On peut itérer $f$
dans le $\C$-espace vectoriel $E_1(\tgamma)$; il vient, pour $a
\in \mbH(\mu)$ et $x \in E_1(\tgamma)$:
$$
f(a\centerdot x) = f(\Int(w_0)(a).x) =
\Int(w_0)(a)\centerdot f(x) = \Int(w_0)^2(a).f(x)=a.f(x).
$$
Donc $f : \twisted{\Int(w_0)}{E_1(\tgamma)} \to
E_1(\tgamma)$ est $\mbH(\mu)$-linéaire et il en découle que
$f^2$ est un automorphisme de $E_1(\tgamma)$. Par le lemme de
Schur $f^2 = \lambda.\id$ pour un $\lambda \in \C^*$. Soit $z
\in \C^*$ tel que $z^2 = \lambda$; quitte à remplacer $f$ par
$\frac{1}{z} f$ on peut ainsi supposer que $f^2 = \id$.
Définissons maintenant
$$
V = \{ x + t_{w_0}.f(x) : x \in E_1(\tgamma)\}.
$$
Soit $x \in E_1(\tgamma)$. Pour $a \in \mbH(\mu)$ il vient
$a.(x + t_{w_0}.f(x)) = a.x + t_{w_0}\Int(w_0)(a).f(x) = a.x +
t_{w_0} f(a.x) \in V$. De plus $t_{w_0}.(x + t_{w_0}.f(x)) =
t_{w_0}.x + f(x) = t_{w_0}.f^2(x) + f(x) \in V$. Donc $V$ est
un sous-$\hd(\bmu)$-module de $F_1(\bgamma)$. Comme
l'application $x + t_{w_0}.f(x) \mapsto x$ de $V$ sur
$E_1(\tgamma)$ est un isomorphisme de $\mbH(\mu)$-modules, $V$
est un sous-module propre de $F_1(\bgamma)$.
\end{dem}

\subsection{Crit\`ere d'irr\'eductibilit\'e de $N(\bgamma)$}  
\label{sec4.5}
On retourne au cas général: $\bgamma = \gamma \otimes \bmu \in
\bar{\calA}$, $\tgamma = \gamma \otimes \bmu \in \calA$. On a
fixé $\sigma \in W_A$ comme au \S~\ref{sec3.1} tel que
$\twisted{\sigma}{\tgamma} = \twisted{\sigma}{\gamma} \otimes
\twisted{\sigma}{\mu}$ avec $\twisted{\sigma}{\mu} = \mu_i$
pour un $i \in \{1,\dots,n\}$.  Rappelons que $R_i$ est le
sytème de racines de type $A_{i-1} \times A_{n-i-1}$ associé,
de groupe de Weyl $W_i = W(R_i) \simeq \mathfrak{S}_i \times
\mathfrak{S}_{n-i}$.  On sait alors que la condition
$\calP_i(\twisted{\sigma}{\gamma}) = \{\alpha \in R_i^+ :
\gamma(\sigma(\alpha)) = \pm 2k \} = \emptyset$ est nécessaire
et suffisante pour que $M(\tgamma)$ soit irréductible, ou
encore que le $\H_{\mathit{gr}}(2k,S_i)$-module
$M(\twisted{\sigma}{\gamma})$ soit irréductible,
cf.~Théorème~\ref{thm13}. Lorsque $n=2m$ on a noté $\tau \in
\GL(\VC)$ l'involution telle que $\tau(\epsilon_j) =
\epsilon_{m+j}$, $\tau(\epsilon_{m+j}) =\epsilon_j$ pour tout
$j\in\{1,\dots, m\}$.

Le théorème suivant donne un critère d'irréductibilité pour
$N(\bgamma)$.

\begin{Theoreme}
\label{thm48}
Avec les notations précédentes:
\begin{enumerate}
\item[{\rm (a)}] si $W_A(\mu) = W_A(\bmu)$,
  $$
  N(\bgamma) \ \text{simple} \ \iff \ 
  \calP_i(\twisted{\sigma}{\gamma})=\emptyset;
  $$
\item[{\rm (b)}] si $W_A(\mu) \neq W_A(\bmu)$, alors $n=2m$,
  $i=m$, et
  $$
  N(\bgamma) \ \text{simple} \ \iff \ 
  \{\text{$\calP_m(\twisted{\sigma}{\gamma})=\emptyset$ et
    $\twisted{\tau}{\gamma} \notin
    W_m.\twisted{\sigma}{\gamma}$}\}.
  $$
\end{enumerate}
\end{Theoreme}

\begin{dem}
  Comme nous l'avons vu au \S~\ref{sec4.4}, l'irréductibilité
  de $N(\bgamma)$ est équivalente à celle du
  $\hd(\bmu_i)$-module $F_1(\twisted{\sigma}{\bgamma})$.
  
  (a) résulte trivialement de~(\ref{eq6}).
  
  (b) Supposons $N(\bgamma)$ simple. Le fait que $M(\tgamma) =
  N(\bgamma)$ comme $\hd$-module force la simplicité de
  $M(\tgamma)$ comme $\hb$-module,
  i.e.~$M(\twisted{\sigma}{\gamma})$ est un
  $\H_{\mathit{gr}}(2k,S_m)$-module simple. Par
  l'Assertion~\ref{ass47} on obtient que
  $M(\twisted{\sigma}{\gamma})$ n'est pas isomorphe à
  $M(\twisted{\tau}{\gamma})$; la Remarque~\ref{rem63} assure
  alors que $\twisted{\tau}{\gamma} \notin
  W_m.\twisted{\sigma}{\gamma}$.
  
  La réciproque est similaire: si $M(\twisted{\sigma}{\gamma})$
  est simple et $\twisted{\tau}{\gamma} \notin
  W_m.\twisted{\sigma}{\gamma}$, la Remarque~\ref{rem63} et
  l'Assertion~\ref{ass47} montrent que
  $F_1(\twisted{\sigma}{\bgamma})$ est irréductible.
\end{dem}

\bigskip

} 
\renewcommand{\refname}{R\'ef\'erences}

\end{document}